\documentclass[11pt]{article}  
\usepackage{geometry}
\usepackage{setspace}
\usepackage{lineno}
\usepackage{graphicx}

\usepackage{amsmath, amssymb, amsthm}

\usepackage{authblk}       
\usepackage{times}  
\usepackage{lmodern}
\usepackage{dcolumn}
\usepackage{bm}
\usepackage[utf8]{inputenc}
\usepackage[T1]{fontenc}
\usepackage{mathptmx}
\usepackage{etoolbox}
\usepackage{subcaption}
\usepackage{array}
\usepackage{float}
\usepackage{hyperref}
\usepackage{multicol}
\usepackage[table]{xcolor}
\usepackage{booktabs}
\usepackage{colortbl}
\usepackage{natbib}
\usepackage{multirow}
\begin{document}
\title{A physics-informed inverse modeling framework for Moose–Wolf dynamics from limited \& noisy data in Isle Royale National Park}

\author[1]{Anurag Singh}
\author[1]{Nitu Kumari}
\affil[1]{School of Mathematical \& Statistical Sciences, IIT Mandi, Himachal Pradesh, 175005, India}

\date{}	
\maketitle 

\section*{Abstract}
The interaction of moose (Alces alces) and wolf (Canis lupus) populations in ecosystems such as Isle Royale National Park is a canonical benchmark for ecological modeling of prey-predator dynamics. Mathematical modeling is a useful tool for modeling these interactions. A central challenge in this domain is the inverse problem, recovering governing system parameters from observational data. Although classical parameter estimation methods have seen considerable progress, they mostly rely on data rather than physics and therefore largely fall short in capturing the time-varying nature of ecological interactions driven by environmental fluctuations, seasonal forcing, and habitat change. This study addresses that gap by solving the inverse problem for a non-autonomous prey-predator system that incorporates $\theta$-logistic prey growth with temporally varying intrinsic growth and natural death rates, as well as Holling type-II and ratio-dependent functional responses. A deep learning framework (self-adaptive bc-PINN with transfer learning) is employed to estimate time-dependent and constant parameters directly from population time-series data (1959--2019). Before estimating the parameters, we performed a structural identifiability analysis to ensure that the model parameters are identifiable. The framework demonstrates good reconstruction and prediction results across both functional responses. The ratio-dependent model has shown the better prediction trend beyond the training data. Our framework successfully predicts the sudden decline in the moose population in $2020$.

\section{Introduction} \label{sec-1}
One of the most important questions in ecology is to understand how the populations of prey and predators influence each other, and this has real-world applications in wildlife management and conservation. Since the pioneering work of \citet{lotka1925elements} and \citet{volterra1926fluctuations}, mathematical models have been central tools for explaining the complex patterns that emerge from these interactions in space and time. One of the best studied examples of such a system is the relationship between moose and wolves in Isle Royale National Park, which has been observed continuously for over six decades, providing an exceptionally rich record of population changes, environmental disturbances, and ecosystem behavior.

Isle Royale, a 544 km² wilderness island in Lake Superior, offers a remarkably isolated and simplified ecosystem that approximates the controlled conditions often assumed in theoretical models. Moose colonized the island in the early twentieth century, followed by wolves in the late 1940s, establishing what \citet{peterson2014trophic} describe as a natural laboratory for studying prey-predator interactions. The winter diet of Isle Royale wolves consists of approximately 90\% moose, and moose have no other significant predators, creating what \citet{jost2005wolves} characterize as a relatively pure two-species prey-predator system. This simplicity, rare among wolf populations that typically inhabit multipredator and multiprey systems, makes Isle Royale uniquely valuable for testing theoretical predictions against empirical observations. The Isle Royale moose-wolf project, initiated in $1958$ and currently led by researchers from Michigan Technological University, has documented remarkable population fluctuations.

A central question in prey-predator theory concerns the appropriate mathematical formulation of the functional response, which quantifies how predation rate changes with prey and predator abundance. Traditional prey-dependent model, exemplified by Holling type-II functional response \citet{holling1959some}, assume that predation rate depends primarily on prey density. However, extensive analysis of Isle Royale data has provided strong empirical support for ratio-dependent functional response \citet{arditi1989coupling}. \citet{vucetich2002effect} found that ratio-dependent models substantially outperformed prey-dependent alternatives in describing wolf predation patterns, a finding reinforced by \citet{jost2005wolves}, who demonstrated that the wolves of Isle Royale display scale-invariant satiation and ratio-dependent predation on moose. This ratio-dependence reflects strong predator interference and competition among wolves, fundamentally altering predictions about system stability and the balance between top-down and bottom-up control.

These ecological facts motivate the study of non-autonomous formulations of prey-predator models, where parameters are time-dependent in order to account for changing environmental conditions. These time-varying parameterizations can include documented shifts in habitat quality, climate-driven changes in prey vulnerability, and temporal variation in predator efficiency. However, the mathematical and computational challenges of parameter estimation increase substantially when moving from autonomous to non-autonomous systems, especially when working with real-world data characterized by measurement noise, missing observations, and inherent stochasticity.

Parameter estimation of systems of differential equations from time series data is an inverse problem \citet{tarantola2005inverse} that has been traditionally approached by methods such as least squares optimization \citet{bard1974nonlinear}, maximum likelihood estimation \citet{lee2025ecological}, and Bayesian inference \citet{coelho2011bayesian}. These approaches also have particular problems, such as ill-conditioned problems, in which different combinations of parameters can give rise to similar dynamics; problems of parameter identifiability, in which the data available are not sufficient to identify all of the parameters; and problems of measurement noise and model misspecification. These problems are particularly acute in ecological systems where observational data are sparse relative to the complexity of the model. However, there are very few system/parameter identification studies on the data available on Isle Royale National Park \citet{knadler2008models, singh2025modeling, lee2025ecological}.

Recent advances in scientific machine learning have led to the emergence of physics-informed neural network (PINN) as a viable alternative paradigm for solving inverse problems in systems governed by differential equations. PINN, introduced by \citet{raissi2019physics}, utilize automatic differentiation to incorporate equation residuals in addition to data-fitting objectives, thereby encoding the physical laws directly into the neural network training process. This method uses physical information as a regularization mechanism and may allow one to accurately estimate parameters even in the case of sparse or noisy observations. The principal novelty is to formulate the inverse problem as a constrained optimization problem where the neural network is required to identify the parameters such that the governing equations are satisfied while simultaneously approximating the solution trajectories. Several researchers have employed this method for both forward and inverse prey-predator problems \citet{slavova2024physics, panchal2025predator, chrisnanto2025unified}, but the std-PINN framework exhibits spectral bias, stiffness, and poor extrapolation issues. Seeing the data of Isle Royale and the nature of the considered equation, the backward compatible physics-informed neural network (bc-PINN) method is chosen to solve our problem.

In this study, we compare the predictive performance of the two prey-predator models: $(i)$ model with Holling type-II functional response and $(ii)$ model with ratio-dependent functional response. These functional responses were chosen based on previous studies of functional responses by \citet{vucetich2002effect, jost2005wolves}, where researchers selected several types of functional responses involved in a parameter estimation analysis from the Isle Royale moose–wolf data. They found that the Holling type-II and ratio-dependent performed best among the considered functional responses at different scales. Also, in place of simple logistic growth, we have taken $\theta$-logistic growth, which is a better approximation in case of large mammals. To estimate the model parameters, we employ a backward-compatible Physics-Informed Neural Network (bc-PINN) \citet{mattey2022novel, kiyani2025optimizing} combined with a self-adaptive weighting strategy \citet{mcclenny2023self} and transfer learning. The observed data exhibit significant fluctuations, and the governing differential equations can become stiff during simulation; these challenges motivate the use of the bc-PINN framework for improved stability and reliability. In addition, the total loss function consists of three distinct components. Achieving an optimal solution requires carefully balancing the contribution of each term, which is effectively handled through the self-adaptive weighting approach that dynamically adjusts the weights during training. 

The main contributions of this study are as follows:
\begin{itemize}
	
	\item We formulate non-autonomous prey--predator models incorporating Holling type-II and ratio-dependent functional responses under a generalized $\theta$-logistic growth framework, enabling a more realistic representation of ecological interactions and environmental variability in Isle Royale National Park.
	
	\item We perform parameter estimation for the proposed models using 61 years of observed time-series data, providing a data-driven assessment of the underlying ecological dynamics.
	
	\item We employ an advanced PINN framework (self-adaptive bc-PINN with transfer learning) capable of simultaneously estimating both time-dependent and constant model parameters from noisy and irregular data, thereby addressing complex inverse problems in ecological modeling.

\end{itemize}

The results show the efficiency of the proposed framework for the identification of the parameters and dynamics of both the prey--predator models. The reconstructed solutions match the observed data very well and demonstrate the power of the framework for model reconstruction and prediction. In addition, the comparison between the two models indicates that the ratio-dependent model provides better prediction trend than the Holling type-II model. The results indicate that the proposed PINN framework provides a reliable and efficient approach for parameter estimation and forecasting in complex ecological systems.

\section{Methodology} \label{sec-2}
Physics-Informed neural network provide a mesh-free method for solving forward and inverse problems. PINN go one step further than traditional neural networks that mainly fit data by incorporating physical laws directly into the model. Inverse problems arise in many areas like engineering, medical imaging, geophysics, epidemiology, and ecology, where the goal is to determine hidden properties of a system using observations. Conventional approaches typically depend on iterative optimization or sophisticated data assimilation techniques, which can be computationally demanding. In contrast, PINN provide a more adaptable and scalable alternative by using deep learning models that naturally embed the governing physical laws of the system into the learning process. PINN framework presents a versatile approach to solve  inverse problems in differential equations. In this work, we focus on the inverse problem of parameter identification for ordinary differential equations (ODEs). The main goal is to solve an inverse problem, which can be difficult and may not have a single clear solution, especially when only limited, noisy, and irregular data is available. 

Consider a general first-order ODE with unknown parameters:
\begin{equation}
\frac{dv}{dt} = f(t, v; \theta), \quad t \in (0, T],
\label{eq-1}
\end{equation}

where $v$ represents the state variable, $t$ is time, $f$ describes the dynamic limitations that establish the equations of motion of the system, $\theta$ denotes the unknown parameters (it can be constant/time dependent), and $T$ is the final time. The objective is to determine $\theta$ given observations of $v$ over time.

For the inverse problem using standard PINN, we define trainable parameters $\theta$ within the neural network framework. The neural network learns both the solution $v(t)$ and the parameter(s) $\theta$ simultaneously by minimizing a composite loss function.

The first component of the loss function, $\mathcal{L}_d(\phi)$ measures the discrepancy between the neural network prediction, $\hat{v}$ and the observed data, $v^o$:
\begin{equation}
\mathcal{L}_d(\phi) = \frac{1}{N_d} \sum_{i=1}^{N_d} \left( \hat{v}(t_i^o) - v_i^o \right)^2, \quad t_i^o \in [0, T],
\label{eq-2}
\end{equation}

where $\phi$ are the parameters (weights and bias) of neural network, $N_d$ is the number of observed data points, $\hat{v}(t_i^o)$ is the neural network output at time $t_i^o$, and $v_i^o$ is the observed value at $t_i^o$. The superscript $({\cdot})^o$ denotes observed data.

The second component, $\mathcal{L}_r(\phi, \theta)$ enforces the governing ODEs by minimizing the residual at the collocation points:

\begin{equation}
\mathcal{L}_r(\phi, \theta) = \frac{1}{N_r} \sum_{k=1}^{N_r} \left(\frac{dv(t_k)}{dt_k} - f(t_k, v(t_k), \theta)\right)^2 , \quad t_k \in [0, T].
\label{eq-3}
\end{equation}

Here $N_r$ is the number of collocation points on which physics residual will be calculated. This term ensures that the learned solution satisfies the physics, which governs the data.

The complete loss function for the inverse std-PINN combines both components:
\begin{equation}
\mathcal{L}(\phi, \theta) = \lambda_d \mathcal{L}_d(\phi) + \lambda_r \mathcal{L}_r(\phi, \theta),
\label{eq-4}
\end{equation}
where, $\lambda_d$ and $\lambda_r$ are hyper parameters. Using an appropriate optimizer, the neural network learns both the model parameters and the solution that best fits the data while adhering to the governing equations.

\subsection{Self-adaptive bc-PINN}
The bc-PINN approach extends the standard PINN framework by solving the inverse problem sequentially over time segments. This strategy offers improved accuracy and efficiency, particularly for long-time integration and complex parameter evolution. The temporal domain $[0, T]$ is partitioned into $n$ segments:

\begin{equation}
[T_0 = 0, ~T_1], ~[T_1, ~T_2], \cdots, ~[T_{n-1}, ~T_n=T],
\label{eq:time_segments}
\nonumber
\end{equation}

where the $i^{th}$ segment is denoted as $\Delta T_i = [T_{i-1}, ~T_i]$ for $i = 1, \ldots, n$.
To estimate the parameters within a given time segment, the residual is minimized using only the observed data from that specific interval. To ensure backward compatibility, the model simultaneously enforces consistency with the solutions obtained from all preceding segments. Furthermore, the parameters for the $n^{th}$ segment are initialized using those learned in the $(n-1)^{th}$ segment, promoting continuity and faster convergence. This approach enables accurate parameter estimation within each segment while also capturing the solution behavior across the entire domain. The corresponding loss function for parameter identification within the self-adaptive bc-PINN framework is defined as follows:
\begin{itemize}
\item Data loss:
\begin{equation}
\mathcal{L}_d(\phi, \lambda_d) = \frac{1}{N_d} \sum_{k=1}^{N_d} [\lambda_d \left( \hat{v}(t_k^o) - v_k^o \right)]^2, \quad t_k^o \in [T_{n-1}, T_n].
\label{eq-5}
\end{equation}

\item Residual loss:
\begin{equation}
\text{Define}~~ \mathcal{R} := \frac{dv}{dt} - f(t, v, \theta).
\nonumber
\end{equation}
\begin{equation}
\mathcal{L}_r(\phi, \theta, \lambda_r) = \frac{1}{N_r} \sum_{k=1}^{N_r} \left[\lambda_r  \mathcal{R}(t_k^r) \right]^2, \quad t_k^r \in [T_{n-1}, T_n].
\label{eq-6}
\end{equation}

\item Backward compatible loss:
\begin{equation}
\mathcal{L}_{bc}(\phi, \lambda_{bc}) = \frac{1}{N_{bc}} \sum_{k=1}^{N_{bc}} \left[\lambda_{bc} \left( \hat{v}(t_k^{bc}) - v_k^{bc} \right)\right]^2, \quad t_k^{bc} \in [0, T_{n-1}],
\label{eq-7}
\end{equation}
where $N_{bc}$ is the number of backward compatible points, $\hat{v}(t_k^{bc})$ is the neural network output, and $v_k^{bc}$ is the observed data at $t_k^{bc}$ . The superscript $({\cdot})^{bc}$ denotes observed data in the previous time segments.

\item Total loss function for inverse self-adaptive bc-PINN is given as
\begin{equation}
\mathcal{L}(\phi, \theta, \lambda_d, \lambda_r, \lambda_{bc}) = \mathcal{L}_d(\phi, \lambda_d) + \mathcal{L}_r(\phi, \theta, \lambda_r) + \mathcal{L}_{bc}(\phi, \lambda_{bc})
\label{eq-8}
\end{equation}
\end{itemize}

\begin{figure}
	\centering
	\includegraphics[width=\textwidth]{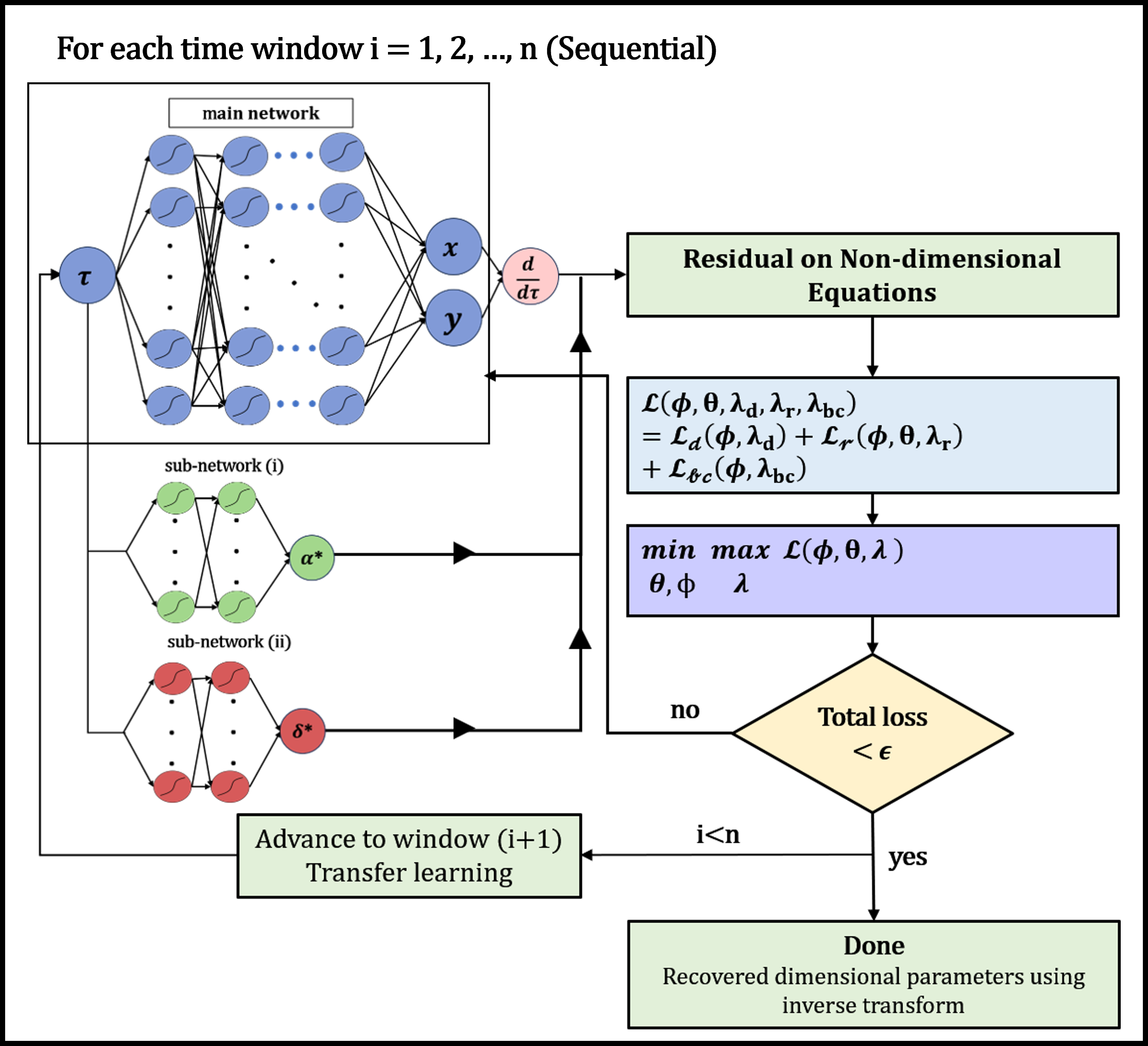}
	\caption{\textbf{Schematic workflow of the self-adaptive bc-PINN framework.} The main network predicts the non-dimensionalize states ($x$, $y$) from the temporal input ($\tau$), and two sub-networks that estimate the parameters $\alpha^*(\tau)$ and $\delta^*(\tau)$. The network outputs are then automatically differentiated to obtain the physics residual on the transformed governing equations.The composite loss function $\mathcal{L}$ is the sum of the data loss ($\mathcal{L}_d$), physics loss ($\mathcal{L}_r$), and backward compatible loss ($\mathcal{L}_{bc}$). The training procedure employs a min max optimization strategy. For the first window, we only have the data loss and physics loss, since the bc loss contributes to the total loss function only from the second window onward. After completion of the entire training procedure, the state variables and parameters from all windows are transformed into physical units using the inverse transformation.}
	\label{fig-1}
\end{figure}

Instead of assigning fixed weights to different loss functions, as in Eq. \eqref{eq-4}, the proposed approach follows the idea of self-adaptation in neural networks. An important aspect of self-adaptive PINN is that the loss function, $\mathcal{L}(\phi, \theta, \lambda_d, \lambda_r, \lambda_{bc})$ is minimized with respect to the network parameters $\phi$ and model parameters $\theta$, while at the same time it is maximized with respect to the adaptive weights $\lambda_d$, $\lambda_r$, and $\lambda_{bc}$. The schematic diagram of self-adaptive bc-PINN with transfer learning framework is presented in Fig.~\ref{fig-1}.

\section{Mathematical models} \label{sec-3}
In this section, we explain the mathematical models used to understand complex system behavior of moose-wolf system. The first model describes the interaction between prey and predator populations using $\theta$-logistic growth and a Holling type-II functional response. The second model represents the same interaction but uses a ratio-dependent response instead.

\subsection{Prey-predator model with Holling type-II functional response}
The first mathematical model used in this study is described by a system of ordinary differential equations:

\begin{equation}
	\begin{aligned}
		\dfrac{dm}{dt} &= \alpha(t)\,m \left(1 - \left(\dfrac{m}{K}\right)^4\right) - \dfrac{a m w}{b + m}, \\[6pt]
		\dfrac{dw}{dt} &= \dfrac{c m w}{b + m} - \delta(t)\,w,
	\end{aligned}
	\label{eq-9}
\end{equation}

where $m(t)$ denotes the prey (moose) population and $w(t)$ denotes the predator (wolf) population at time $t$. The prey grows logistically with intrinsic growth rate $\alpha (t)$ and carrying capacity $K$. Predation follows a Holling type-II functional response, where $a$ is the the maximum per-capita predation rate (or maximum feeding rate) and $b$ is the half-saturation constant, the prey density at which $a$ is half its maximum. The predator grows through conversion of consumed prey with efficiency $c$, and dies at a natural rate $\delta (t)$.

To non-dimensionalize the system, we introduce the non-dimensional variables
\begin{equation}
	x=\dfrac{m}{K}, \qquad
	y=\dfrac{w}{W_c}, \text{ and }  \qquad
	\tau=\dfrac{t-t_0}{T_c}.
	\label{eq-10}
\end{equation}

Substituting these transformations into Eq.~\eqref{eq-9} gives the non-dimensional model
\begin{equation}
	\begin{aligned}
		\dfrac{dx}{d\tau} &= \alpha^*(\tau) x  (1-x^4) - \dfrac{\tilde a\,x y}{\tilde b + x}, \\
		\dfrac{dy}{d\tau} &= \dfrac{\tilde c\,xy}{\tilde b + x} - \delta^*(\tau)y.
	\end{aligned}
	\label{eq-11}
\end{equation}

The dimensionless parameters are
\begin{equation}
	\alpha^*(\tau) = T_c \alpha(\tau T_c + t_0), \quad
	\tilde a=\dfrac{T_c W_c a}{K},
	\quad
	\tilde b=\dfrac{b}{K},
	\quad
	\tilde c = c T_c, \text{ and } \quad
	\delta^*(\tau) = T_c \delta(\tau T_c + t_0).
	\label{eq-12}
\end{equation}

Using the inverse relations, the corresponding dimensional parameters are given by

\begin{equation}
	m = K x, \quad w = W_c y, \quad  t = \tau T_c + t_0,\quad
	\alpha(t)=\frac{\alpha^*(\tau)}{T_c}, \quad
	a=\frac{\tilde a\,K}{T_cW_c},
	\quad
	b=\tilde b\,K,
	\quad
	c=\frac{\tilde c}{T_c}, \text{ and } \quad
	\delta(t)=\frac{\delta^*(\tau)}{T_c}.
	\label{eq-13}
\end{equation}

\subsection{Model with ratio-dependent functional response}
The second mathematical model is given as:
\begin{equation}
	\begin{aligned}
		\dfrac{dm}{dt} &= \alpha(t)\,m \left(1 - \left(\dfrac{m}{K}\right)^4\right) - \dfrac{a m w}{b w + m}, \\[6pt]
		\dfrac{dw}{dt} &= \dfrac{c m w}{b w + m} - \delta(t)\,w.
	\end{aligned}
	\label{eq-14}
\end{equation}
In the ratio-dependent formulation, the predation rate depends on the ratio $m/w$ rather than on prey density alone, capturing the effect of predator interference. The term $\dfrac{am}{bw + m}$ is the ratio-dependent functional response, where $b$ represents the specific ratio of prey to predators ($m/w$) at which the maximum per-capita predation rate is exactly half of its maximum value. All other parameter descriptions are same as in model (\ref{eq-9}).

To non-dimensionalize the system, we introduce the non-dimensional variables
\begin{equation}
	x=\dfrac{m}{K}, \qquad
	y=\dfrac{w}{W_c}, \text{ and }  \qquad
	\tau=\dfrac{t-t_0}{T_c}.
	\label{eq-15}
\end{equation}

Substituting these transformations into Eq. (\ref{eq-14}) gives the non-dimensional model
\begin{equation}
	\begin{aligned}
	\dfrac{dx}{d\tau} &= \alpha^*(\tau) x  (1-x^4) - \dfrac{\tilde a\,x y}{\tilde b\,y + x}, \\
	\dfrac{dy}{d\tau} &= \dfrac{\tilde c\,xy}{\tilde b\,y+x} - \delta^*(\tau)y.
	\end{aligned}
	\label{eq-16}
\end{equation}

The dimensionless parameters are
\begin{equation}
	\alpha^*(\tau) = T_c \alpha(\tau T_c + t_0), \quad
	\tilde a=\dfrac{T_c W_c a}{K},
	\quad
	\tilde b=\dfrac{bW_c}{K},
	\quad
	\tilde c = c T_c, \text{ and } \quad
	\delta^*(\tau) = T_c \delta(\tau T_c + t_0).
	\label{eq-17}
\end{equation}

Using the inverse relations, the corresponding dimensional parameters are given by

\begin{equation}
	m = K x, \quad w = W_c y, \quad  t = \tau T_c + t_0,\quad
	\alpha(t)=\frac{\alpha^*(\tau)}{T_c}, \quad
	a=\frac{\tilde a\,K}{T_cW_c},
	\quad
	b=\frac{\tilde b\,K}{W_c},
	\quad
	c=\frac{\tilde c}{T_c}, \text{ and } \quad
	\delta(t)=\frac{\delta^*(\tau)}{T_c}.
	\label{eq-18}
\end{equation}

\textbf{\textit{Remark:}} These mathematical models have not been chosen arbitrarily; rather, they have been constructed based on the available data and previous studies on the functional responses and growth rates of the Isle Royale ecosystem. Ecological systems are frequently subject to periodic or irregular environmental forcing--such as seasonal fluctuations in resource availability or climate-driven shifts in mortality rates, which render model parameters inherently time-dependent. This behavior is well documented in the Isle Royale moose-wolf system. Accordingly, in both models considered here, the prey intrinsic growth rate $\alpha(t)$ and the predator death rate $\delta(t)$ are treated as time-varying quantities. Furthermore, in place of classical logistic growth, we adopt the $\theta$-logistic growth framework, with $\theta = 4$, which captures the nonlinear nature of density dependence by confining its effect primarily to population levels close to the carrying capacity--a biologically realistic assumption for large mammal systems such as moose \citet{gilpin1973global, vucetich2011predicting}.

\subsection{Structural identifiability}
Determining identifiability is the natural prerequisite toward reliable parameter estimation. In a parametric ODE model, a parameter is said to be structurally globally identifiable if its value can be uniquely determined from the input-output data, assuming noise-free data and sufficiently informative inputs. In this study, we are performing parameter estimation on Eqs.~\eqref{eq-11} and \eqref{eq-16}. Before parameter estimation, we performed a structural identifiability analysis of these models by assuming that $\alpha^*(\tau)=\alpha_0$ and $\delta^*(\tau)=\delta_0$ are scalar parameters rather than time-dependent functions. The analysis was carried out using the Julia code provided by \citet{dong2023differential}. Given observations of the state variables $x$ and $y$, we found that all scalar parameters, $\alpha_0$, $a$, $b$, $c$, and $\delta_0$ are globally identifiable for both models. This result indicates that the model parameters can, in principle, be uniquely determined from ideal noise-free observations of the state variables.

In our models $\alpha^*$ and $\delta^*$ are time-dependent. But, if we choose one particular time, say $\tau_0$, and consider $\alpha^*(\tau_0)$ and $\delta^*(\tau_0)$ as fixed values. Then, the system becomes the same autonomous model that we have already analyzed. Since that model was found to be globally identifiable, we can say that there is no redundancy in the model structure. In simple terms, different values of $a, b, c, \alpha^*(\tau) \text{ and } \delta^*(\tau)$ cannot produce exactly the same $x$ and $y$ solutions when the forcing values are fixed. Therefore, this analysis serves as a necessary preliminary check before proceeding with parameter estimation. However, it does not guarantee the uniqueness of the estimated solution, particularly in the presence of noisy data.

\section{Training setup} \label{sec-4}
In this section, we describe the training procedure used for self-adaptive bc-PINN to solve the inverse problem.

\subsection{Data and non-dimensionalization} 
\begin{figure}[ht!]
    \centering
    \includegraphics[width=\textwidth, trim=0 0 0 6cm, clip]{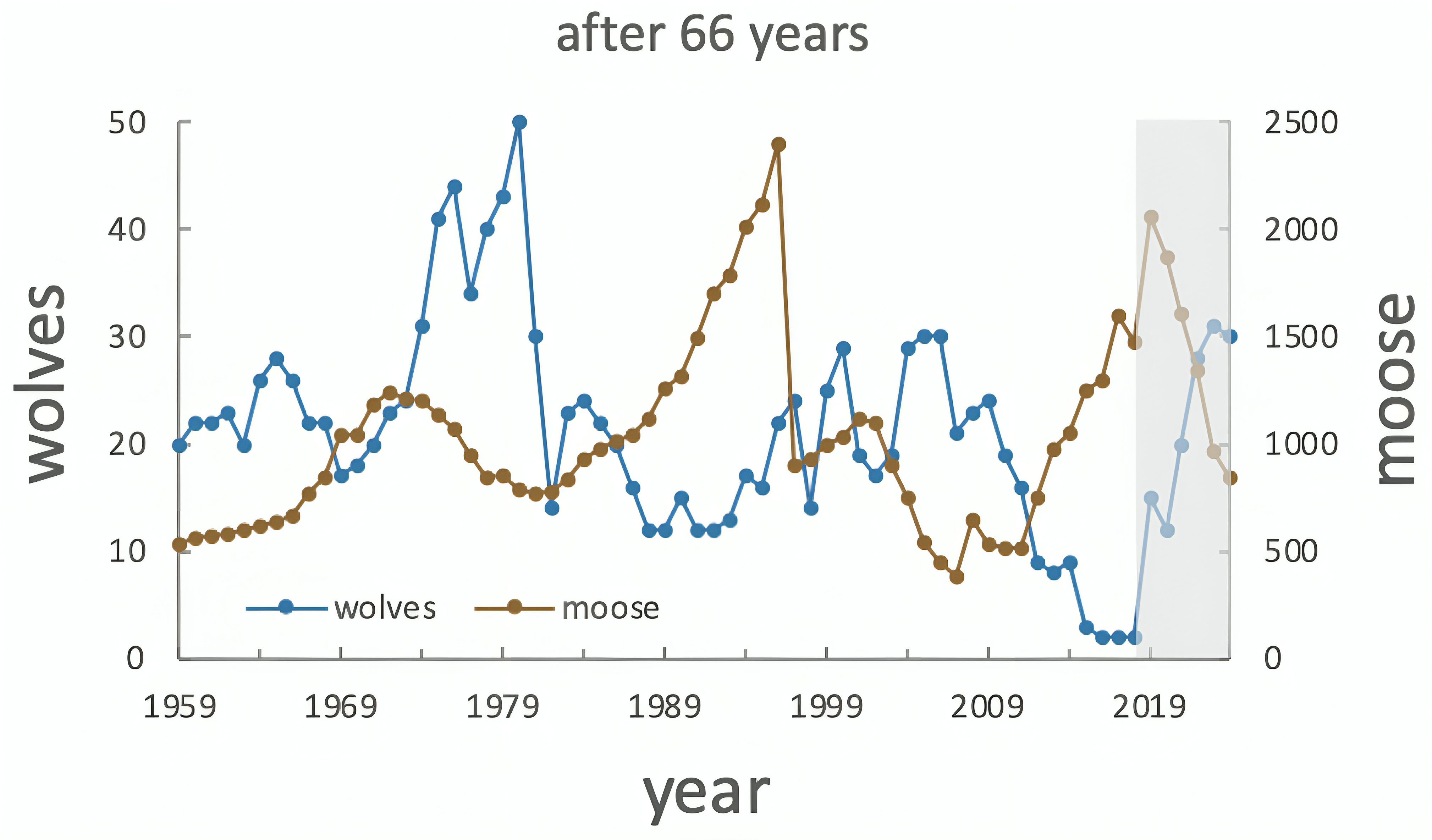}
    \caption{Time series data of moose and wolves in Isle Royale National Park, representing 67 years (1959–2024) of population fluctuations (\href{https://www.isleroyalewolf.org}{https://www.isleroyalewolf.org}). The numerical data is only available from 1959--2019; after 2019, only graphical data is available.}
    \label{fig-2}
\end{figure}

The empirical data used in this work are shown in Fig.~\ref{fig-2}. Over this six-decade period, the wolf population exhibited considerable variability, fluctuating between a minimum of 2 individuals and a maximum of 50, with a mean value of 21.0655 and a standard deviation of 9.9158. In contrast, the moose population density showed a broader range, varying from 385 to 2398, and was characterized by an average population size of 1020.6393 and a standard deviation of 430.9177. Since the data were collected via aerial survey, there may be some uncertainty in the data. Moose and wolves have the same unit: individuals.

\begin{figure}[ht!]
	\centering
	\includegraphics[width=\textwidth, height=10 cm]{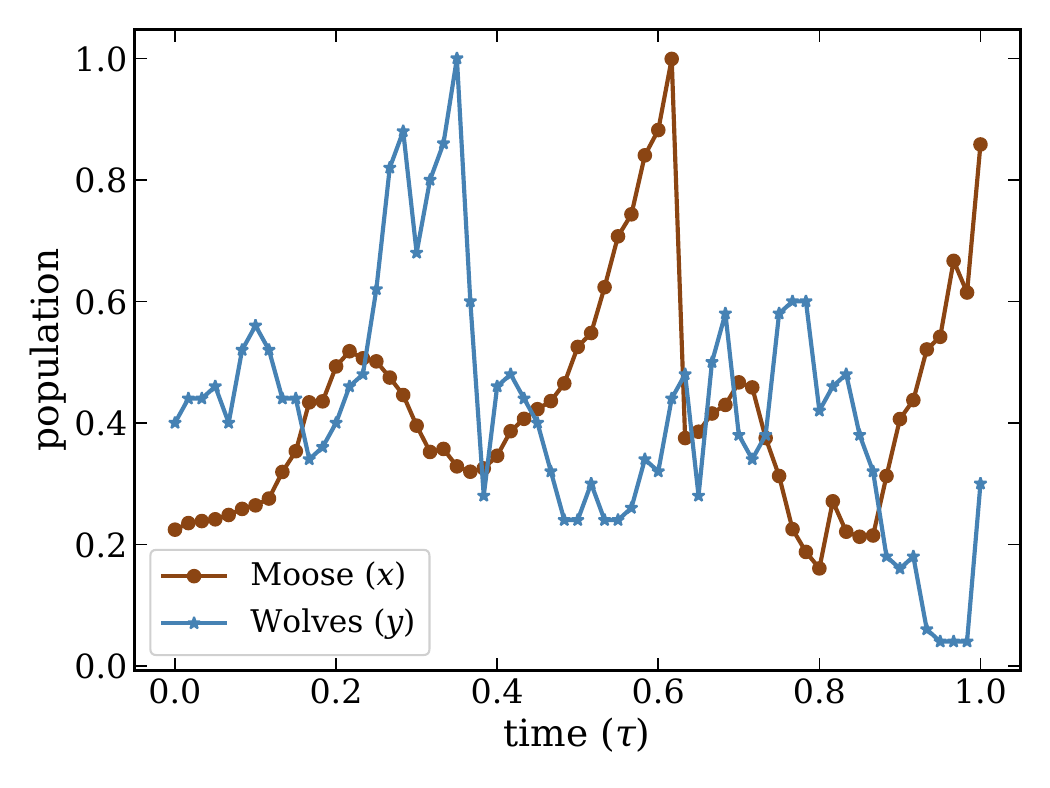}
	\caption{\textbf{Non-dimensionalized data:} The data were obtained by applying the transformations
		$ x = \dfrac{m}{K}, \; y = \dfrac{w}{W_c}, \; \text{and} \; \tau = \dfrac{t - t_0}{T_c} $ to the $61$-year observational dataset with the values $ K = 2400, ~W_c = 50, ~t_0 = 1959, \text{ and } T_c = 60 $.}
	\label{fig-3}
\end{figure}

We have used Eq.~\eqref{eq-10} with the values $K = 2400, ~W_c = 50, ~t_0 = 1959, \text{ and } T_c = 60 $ to non-dimensionalize the data. The carrying capacity and the scaling factor used to non-dimensionalize the wolf population were assumed to be approximately equal to the maximum observed moose and wolf populations, respectively, over the $67$-year study period. Here, $t_0$ denotes the initial year of observation and $T_c$ denotes the total time span of the observations. Now time ($\tau$) and both populations ($x \text{ and } y$) are of order one and lie in $[0, 1] $ (see Fig.~\ref{fig-3}). Training is carried out entirely in these variables, and the recovered quantities are mapped back to physical units afterwards using Eq. \eqref{eq-13} for type-II  and Eq.~\eqref{eq-18} for ratio dependent. The rescaling matters because the raw populations differ by two order of magnitude, and a loss formed from their squared errors would otherwise be governed almost entirely by the moose residual.

\subsection{Segmentation of the domain}
The time interval $[0, 1]$ is partitioned into $15$ segments. Segment $i$ occupies $[T_{i-1}, ~T_i]$ and carries its own set of observations. Three point sets are attached to segment $i$:
\begin{itemize}
	\item Data points: The $N_d$ observations lying inside the segment, at which the network is compared with the real data.
	
	\item Collocation points: $N_r = 150$ points spaced uniformly across the segment, at which governing equations are enforced.
	
	\item Backward compatibility points: The observations of every preceding segment, that is all data points with $\tau \in [0,~ T_{i-1}]$. Their number is $N_{bc} = 4 (i-1)$, which grows from $0$ in the first segment to $56$ in the last.

\end{itemize}  

\subsection{Network architecture}
Three fully connected networks are used. All used the $Tanh$ activation, and all weights are initialized by Xavier initialization.

The state network takes $\tau$ as input and produces the output $(x, y)$. It has $5$ hidden layers with $20$ neurons. The sub-network $(i)$ and $(ii)$ takes $\tau$ and returns a single value: sub-network $(i)$ gives the dimensionless growth rate $\alpha^{*}(\tau)$ and sub-network $(ii)$ gives the dimensionless death rate $\delta^{*}(\tau)$. Both sub-networks have $2$ hidden layers of width $5$. A $softplus$ activation is applied to the output of the sub-network $(ii)$ so that $\delta^{*}$ is positive.

Three scalar coefficients of the model are unknown and are recovered jointly with the network weights. The parameter $a$ is confined to the interval $[a_{\min}, ~a_{\max}] = [18, 21]$ in physical units. So we bound the value of $\tilde{a}$ using the following equation:

\begin{equation}
	\tilde{a} = \tilde{a}_{\min}
	+ \left(\tilde{a}_{\max} - \tilde{a}_{\min} \right)
	\sigma\!\left( a_{\mathrm{raw}} \right),
	\qquad
	\sigma(z) = \frac{1}{1 + e^{-z}},
	\label{eq-19}
\end{equation}

where $\tilde{a}_{\min}$ and $\tilde{a}_{\max}$ denote the interval endpoints in
dimensionless form. Since $\sigma : \mathbb{R} \to (0,1)$, the image of
Eq.~\eqref{eq-19} is contained in $(\tilde{a}_{\min}, ~\tilde{a}_{\max})$ for any
value of $a_{\mathrm{raw}}$. The remaining two coefficients are required only to be
positive and are obtained through the softplus map,
\begin{equation}
\tilde{b} = \ln\!\left( 1 + e^{\,b_{\mathrm{raw}}} \right),
	\qquad
\tilde{c} = \ln\!\left( 1 + e^{\,c_{\mathrm{raw}}} \right),
	\label{eq-20}
\end{equation}
whose range is $(0,\infty)$. The optimizer acts on $a_{\mathrm{raw}}$,
$b_{\mathrm{raw}}$, and $c_{\mathrm{raw}}$ alone, and the constrained coefficients are
recovered from Eqs.~\eqref{eq-19} and \eqref{eq-20} whenever the loss
is evaluated. The loss is minimized over the network weights $\phi$, model coefficients $\theta$, and simultaneously maximized over the point weights $\lambda$,
\begin{equation}
	\min_{\phi,\,\theta}\ \max_{\lambda}\ \mathcal{L}(\phi,\theta,\lambda).
	\label{eq-21}
\end{equation}

Both directions are realized within a single optimizer by a gradient-reversal layer. One Adam step then descends in $(\phi,\theta)$ and ascends in $\lambda$. \textit{Adam} is applied with an initial learning rate of $5\times 10^{-3}$, decayed by a cosine schedule to $10^{-5}$ over the number of Adam iterations assigned to each segment, and without weight decay. The gradient norm of $(\phi,\theta)$ is clipped at unity.

Segments are solved consecutively in order of increasing time. The first is initialized by Xavier normal initialization, with the parameter $\tilde{a}$ at the center of its admissible interval and the remaining two parameters at unity. Each subsequent segment inherits the converged state of its predecessor: the weights of all three networks are copied, and the model coefficients are carried over as starting values. The point weights
alone are reset to unity, since they are attached to point sets that differ between segments. Retention of earlier information is enforced by the backward-compatibility term. Segment $i$ is required to reproduce the observations of all preceding segments, and segment $(i-1)$ was subject to the same condition over its own history. The requirement therefore chains recursively, and the solution remains valid over the entire domain solved up to that point.

\section{Results} \label{sec-5}
This section presents the findings of our study. We begin with our inverse self-adaptive bc-PINN implementation applied to the Holling type-II functional response model, where we present the results and discuss the overall performance. We then shift our attention to the ratio dependent model, following the same evaluative approach. Together, both models are compared in terms of how well they reconstruct the underlying dynamics and how reliably they generalize beyond the training data, offering a clearer picture of their respective strengths and limitations.

\textit{\textbf{System setup:}} All numerical experiments were performed using Python~$3.12.4$ and the PyTorch~$2.9.1$ library with double-precision arithmetic (\texttt{torch.float64}). The computations were carried out on a CPU-based platform running the Windows~$10$ operating system. Our code is also compatible with CUDA; however, due to the computational resources available to us, all computations were performed using the CPU.

\subsection{Model with Holling type-II functional response}
\begin{figure}[ht!]
	\centering	
	\begin{subfigure}{0.49\textwidth}
		\centering
		\includegraphics[width=\linewidth]{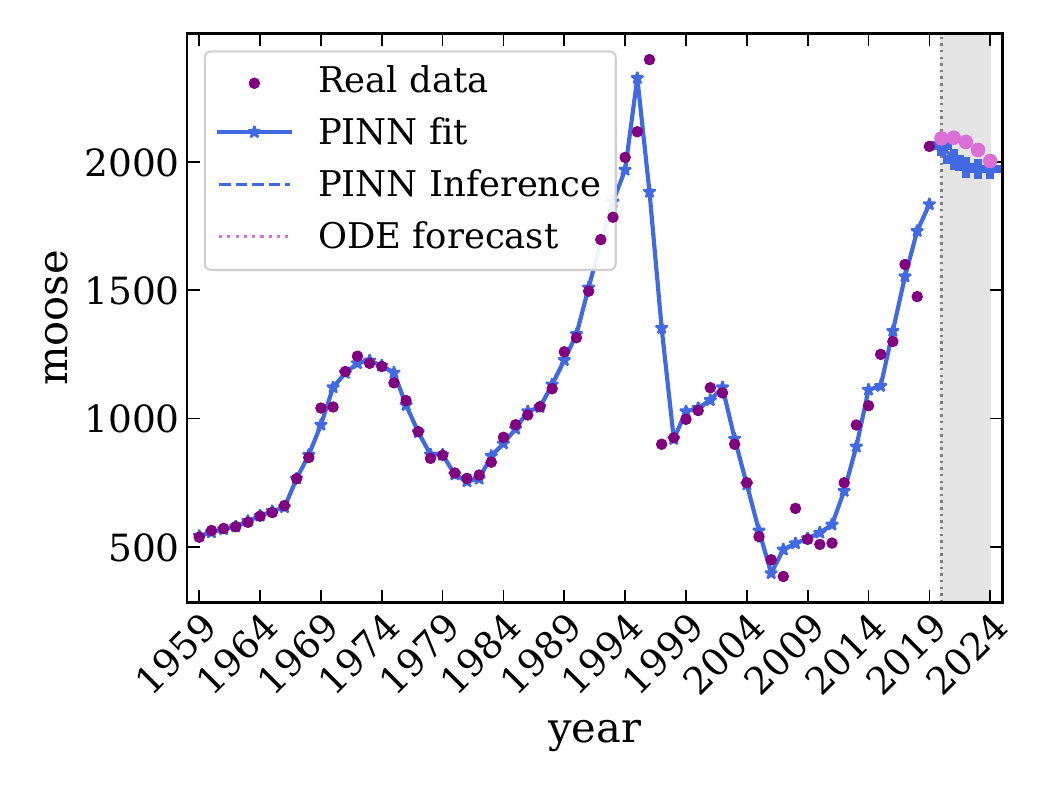}
		\caption{}
		\label{fig-4(a)}
	\end{subfigure}
	\hfill
	\begin{subfigure}{0.49\textwidth}
		\centering
		\includegraphics[width=\linewidth]{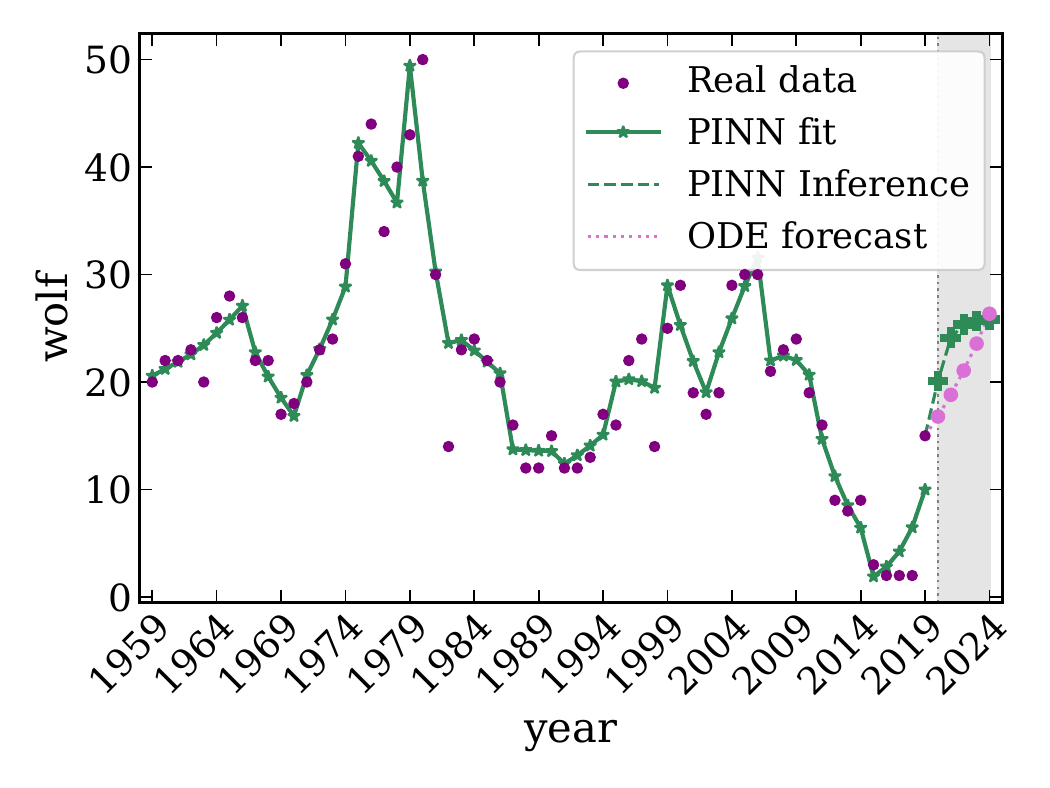}
		\caption{}
		\label{fig-4(b)}
	\end{subfigure}	
	\caption{\textbf{Holling type--II:} Observed and fitted population trends (1959--2019), including predictions for 2020--2024 (gray shaded region). Panel \textbf{(a)} shows the moose population dynamics, while panel \textbf{(b)} shows the wolf population dynamics.}
	\label{fig-4}
\end{figure}

The neural network for the Type-II model was trained using $10000$ Adam iterations for each segment, resulting in a total of $150000$ training iterations.
From Fig.~\ref{fig-4}, it can be observed that the PINN model accurately interpolates the real data and is also capable of forecasting population dynamics beyond the training period. Fig.~\ref{fig-4(a)} presents the predicted moose population for the years 2020--2024, obtained using both the neural network and the ODE model with the parameters estimated from the last training segment. The predicted trajectory captures the overall trend of the observed moose population. Similarly, Fig.~\ref{fig-4(b)} shows that the predicted wolf population follows the general pattern of the observed data during the forecasting period. The actual population trends of moose and wolves for the years 2020--2024 are shown in Fig.~\ref{fig-2}. Although the neural network is able to fit the observed data, the ODE model is preferred for forecasting because it incorporates the underlying ecological interactions governing the system. The ODE based predictions are more interpretable and physically consistent, whereas neural networks may fail to generalize outside the training interval. Consequently, ODE models are expected to provide more reliable long-term forecasts.

\begin{figure}[ht!]
	\centering	
	\begin{subfigure}{0.49\textwidth}
		\centering
		\includegraphics[width=\linewidth]{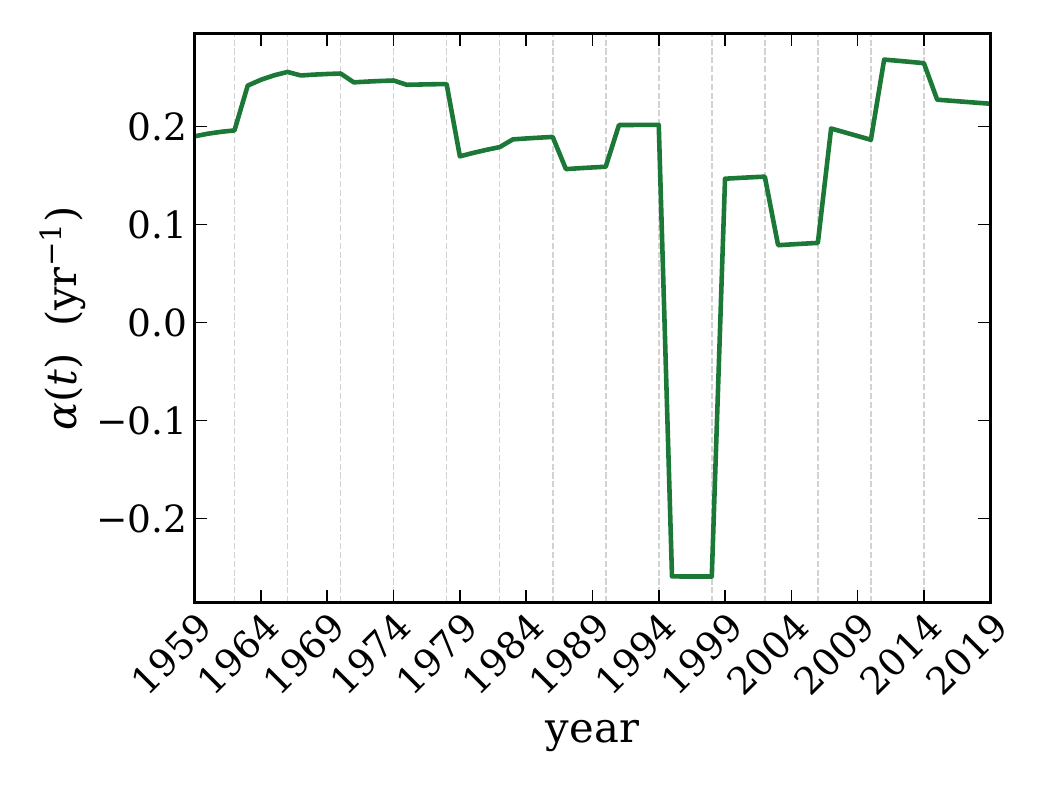}
		\caption{}
		\label{fig-5(a)}
	\end{subfigure}
	\hfill
	\begin{subfigure}{0.49\textwidth}
		\centering
		\includegraphics[width=\linewidth]{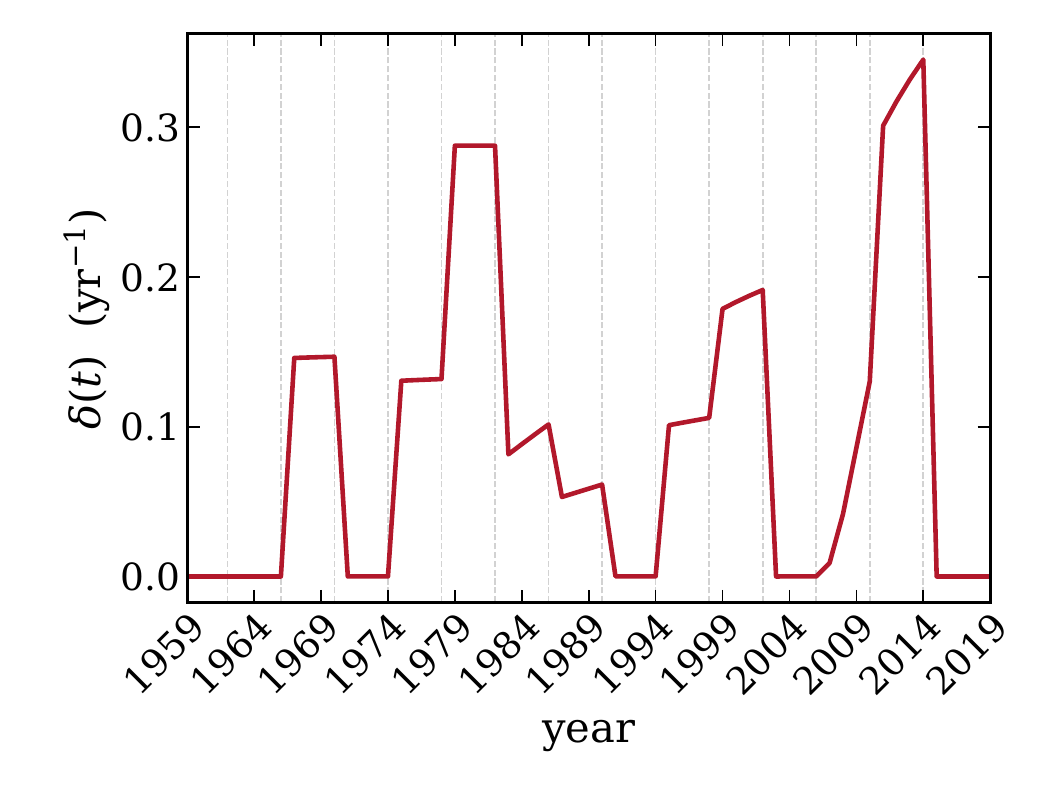}
		\caption{}
		\label{fig-5(b)}
	\end{subfigure}
	\caption{Estimated time dependent parameters: \textbf{(a)} $\alpha(t)$ and \textbf{(b)} $\delta(t)$ of \textbf{Holling type-II model}.}
	\label{fig-5}
\end{figure}

Fig.~\ref{fig-5} shows the estimated time dependent parameters, i.e. intrinsic growth rate of prey using sub-network (i) and death rate of predator using sub-network (ii). Fig.~\ref{fig-5(a)} shows that the maximum growth rate is $\approx 0.2682$ and minimum growth rate is $\approx -0.2592$. The sudden decrease in growth rate in $10^{th}$ window is due to many environmental factors, which also can be seen in Fig.~\ref{fig-2}. From Fig.~\ref{fig-5(b)}, it can be observed that the death rate varies from $10^{-8}$ to $0.3451$. The sharp decline in the wolf population around $1981$ was primarily caused by the outbreak of canine parvovirus; however, this effect has not been incorporated into our model due to the lack of available data. Furthermore, the reintroduction of wolves during 2018--2019 is also not included in the model formulation.

\begin{figure}[ht!]
	\centering	
	\begin{subfigure}{0.32\textwidth}
		\centering
		\includegraphics[width=\linewidth]{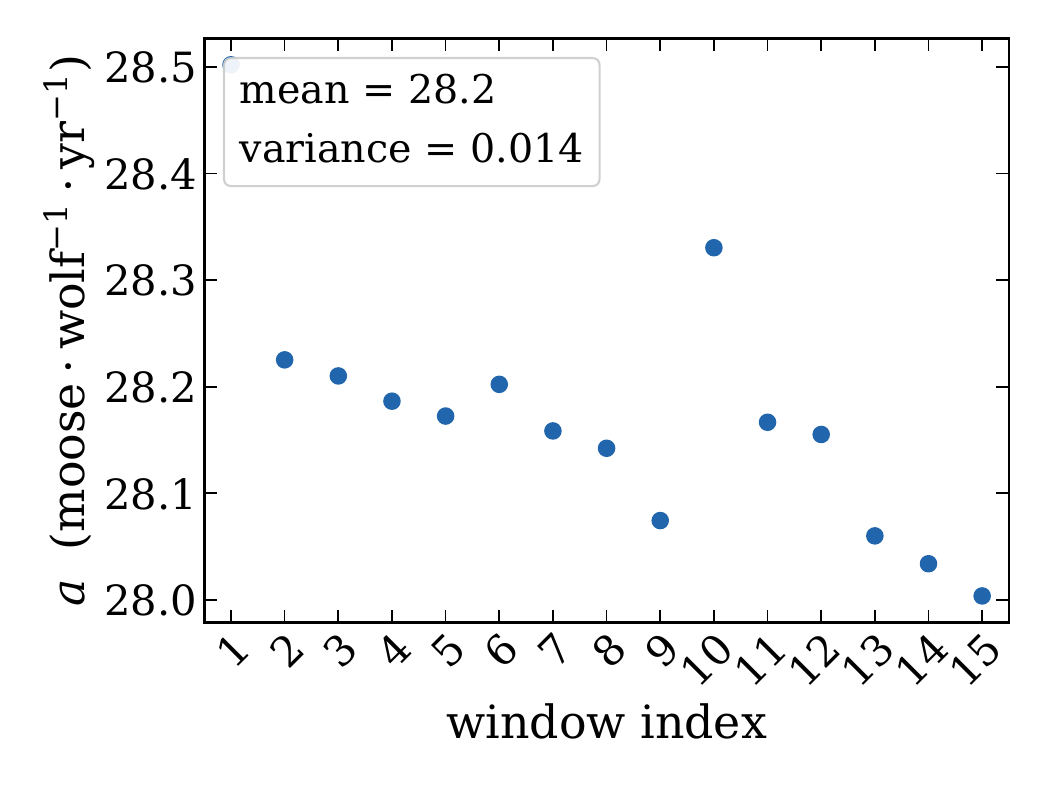}
		\caption{}
		\label{fig-6(a)}
	\end{subfigure}
	\hfill
	\begin{subfigure}{0.32\textwidth}
		\centering
		\includegraphics[width=\linewidth]{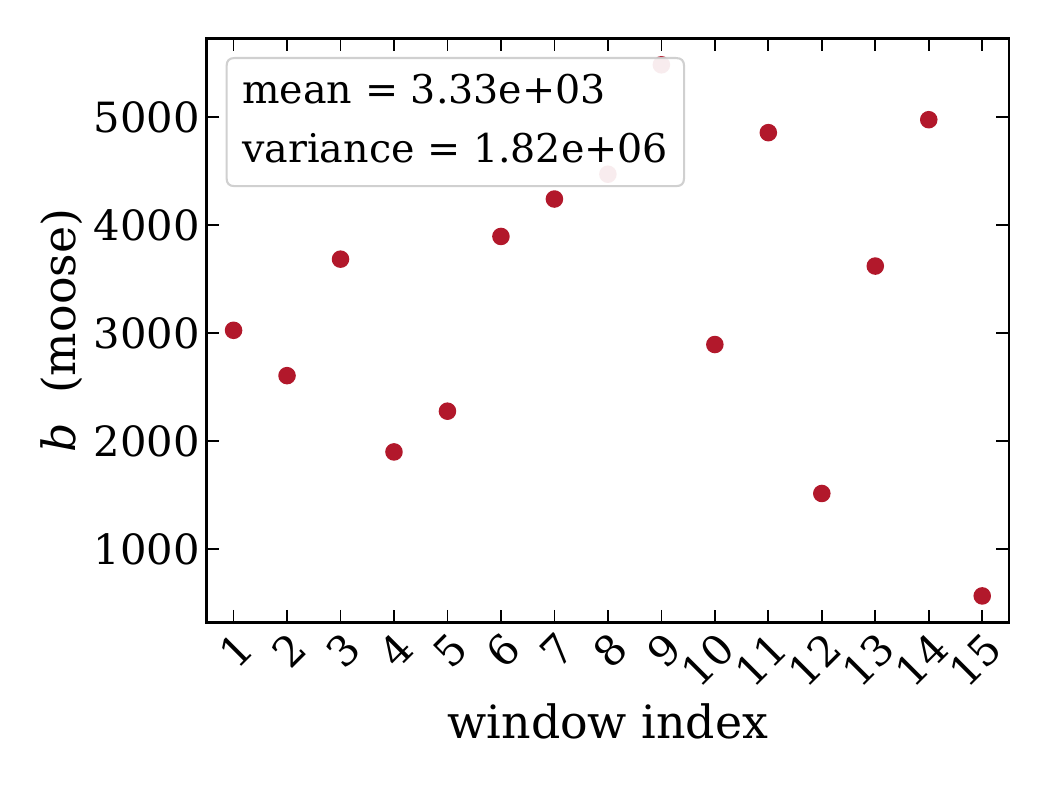}
		\caption{}
		\label{fig-6(b)}
	\end{subfigure}
	\hfill	
	\begin{subfigure}{0.32\textwidth}
		\centering
		\includegraphics[width=\linewidth]{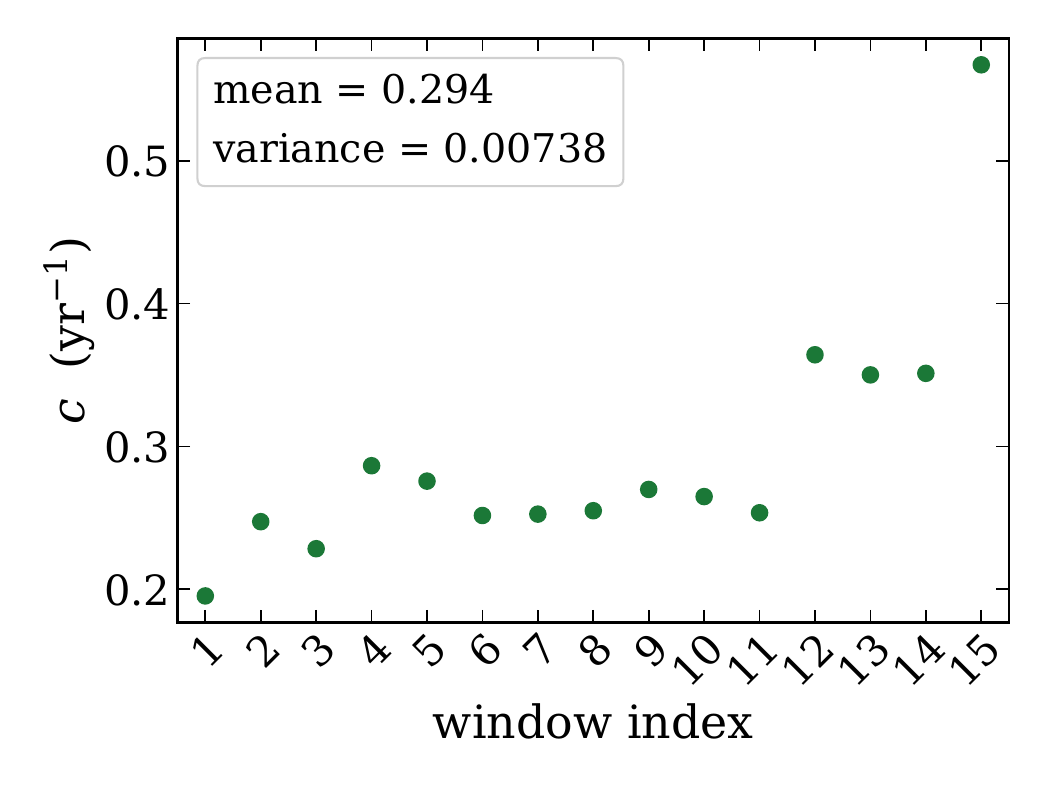}
		\caption{}
		\label{fig-6(c)}
	\end{subfigure}
	\caption{Constant parameters of the \textbf{Holling type-II model} recovered independently in each of
		the $15$ time segments, with the mean and variance across segments inset in each panel.
		\textbf{(a)} The parameter $a$ comes out almost the same in every segment: the estimates span $28.004$ to
		$28.502$ with variance of $0.114$ and mean $28.2$.
		\textbf{(b)} The half-saturation constant $b$, scatters from $562.548$ to $5484.503$ with a standard
		deviation of $1340.55$ and mean $3.33\times10^{3}$. \textbf{(c)} The conversion rate $c$, lies between $0.1953$ and $0.5672$, with a variance of $0.0073$ against a mean of $0.294$, and drifts upward in
		the later segments. The mean and variance is taken over the $15$ segment estimates.}
	\label{fig-6}
\end{figure}

Fig.~\ref{fig-6} shows how tightly the data fix each constant. The maximum per capita
predation rate, $a$ is pinned down almost exactly, every segment returning a value
within about one percent of $28.2$. The conversion rate $c$ rises in the last segment because the sudden increase in the wolf population. This growth in the population of wolves is not due to natural growth but it is due to reintroduction of wolves in the park. The half-saturation constant $b$ is not fixed by the data at
all: its estimates scatter across nearly an order of magnitude with no trend. 

\begin{figure}[ht!]
	\centering	
	\begin{subfigure}{0.49\textwidth}
		\centering
		\includegraphics[width=\linewidth]{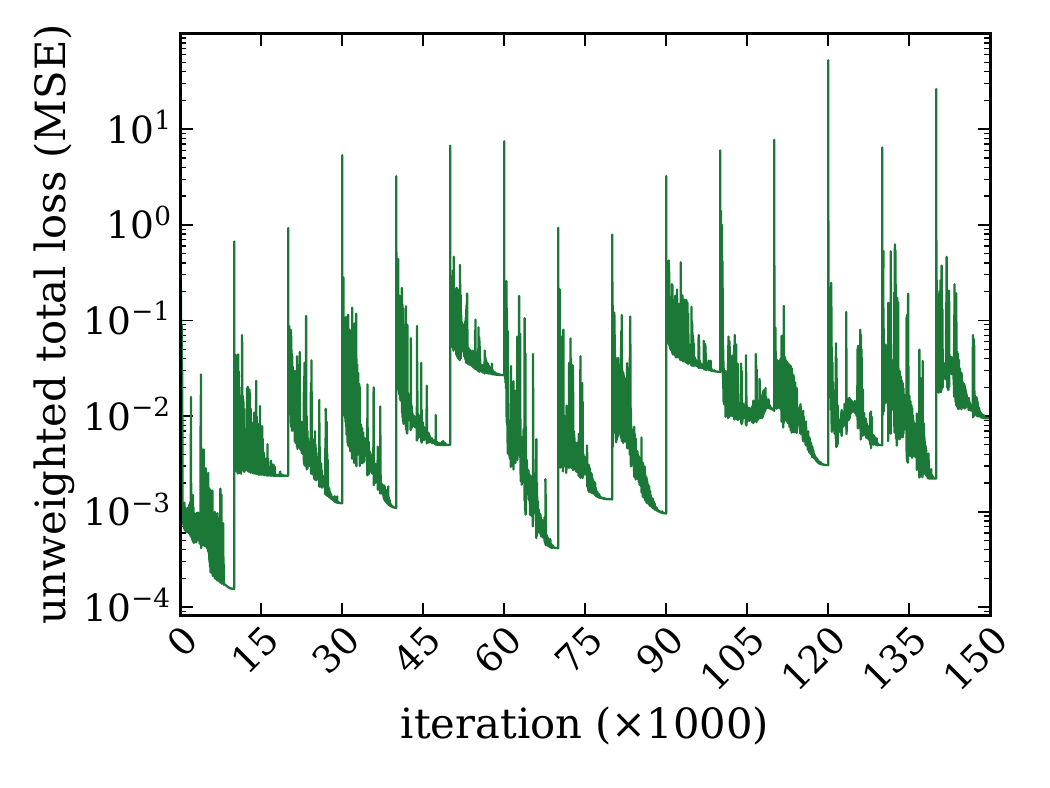}
		\caption{}
		\label{fig-7(a)}
	\end{subfigure}
	\hfill
	\begin{subfigure}{0.49\textwidth}
		\centering
		\includegraphics[width=\linewidth]{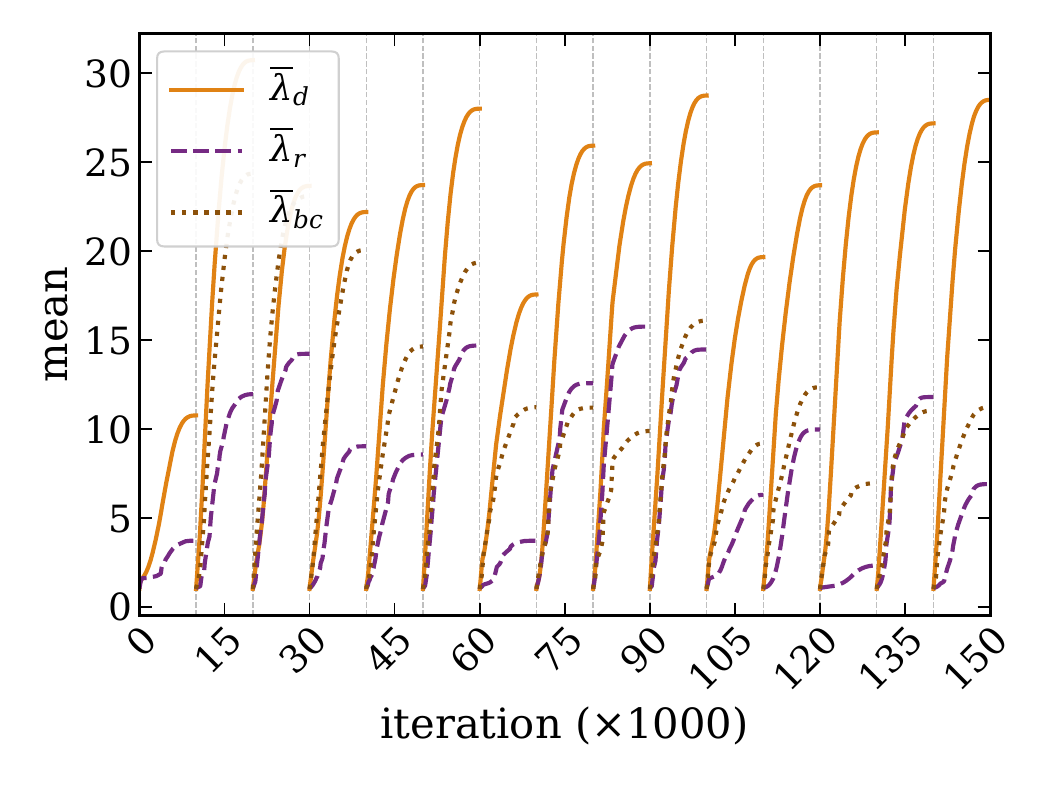}
		\caption{}
		\label{fig-7(b)}
	\end{subfigure}
	\caption{Convergence history of the sequential sweep for \textbf{Holling type-II model}, plotted against the cumulative
		iteration count over all $15$ segments ($10000$ Adam iterations each, $150000$ in
		total).
		\textbf{(a)} Unweighted total loss, $\mathcal{L}^{\mathrm{unw}}=\mathcal{L}^{\mathrm{unw}}_{d}+
		\mathcal{L}^{\mathrm{unw}}_{r}+\mathcal{L}^{\mathrm{unw}}_{bc}$, obtained by setting
		$\lambda\equiv1$, on a logarithmic scale. Within each segment the error falls. The transient spike at each boundary is the cost of moving the
		inherited network onto a new interval, and its decreasing height in later segments shows
		that transfer learning progressively improves the starting point.
		\textbf{(b)} The average weight given to each of the three parts of the loss: data term, physics residual term, and backward-compatible term. These weights are not fixed
		numbers. During training the weight of any point that is still badly fitted is pushed
		upward, so every weight climbs while a segment is being solved and then drops back to one
		when the next segment begins.}
	\label{fig-7}
\end{figure}

Fig.~\ref{fig-7} presents the convergence of the bc-PINN and evolution of the mean weight given to the each loss component. Each point in a loss term carries its own weight, so a segment containing $N$ points has $N$
separate weights rather than a single one. To display them as one curve we take their
arithmetic mean,
\begin{equation}
	\overline{\lambda} = \frac{1}{N}\sum_{k=1}^{N} \lambda_{k},
	\label{eq-22}
\end{equation}
computed separately for the data, physics, and backward-compatible terms and recorded at
every iteration. The mean is used because the three terms contain different numbers of
points, so a sum would scale with $N$ and make the terms incomparable. Every weight is set to
one at the start of a segment, so each curve begins at $\overline{\lambda}=1$ (see Fig.~\ref{fig-7(b)}).

\subsection{Model with ratio-dependent functional response}
\begin{figure}[ht!]
	\centering	
	\begin{subfigure}{0.49\textwidth}
		\centering
		\includegraphics[width=\linewidth]{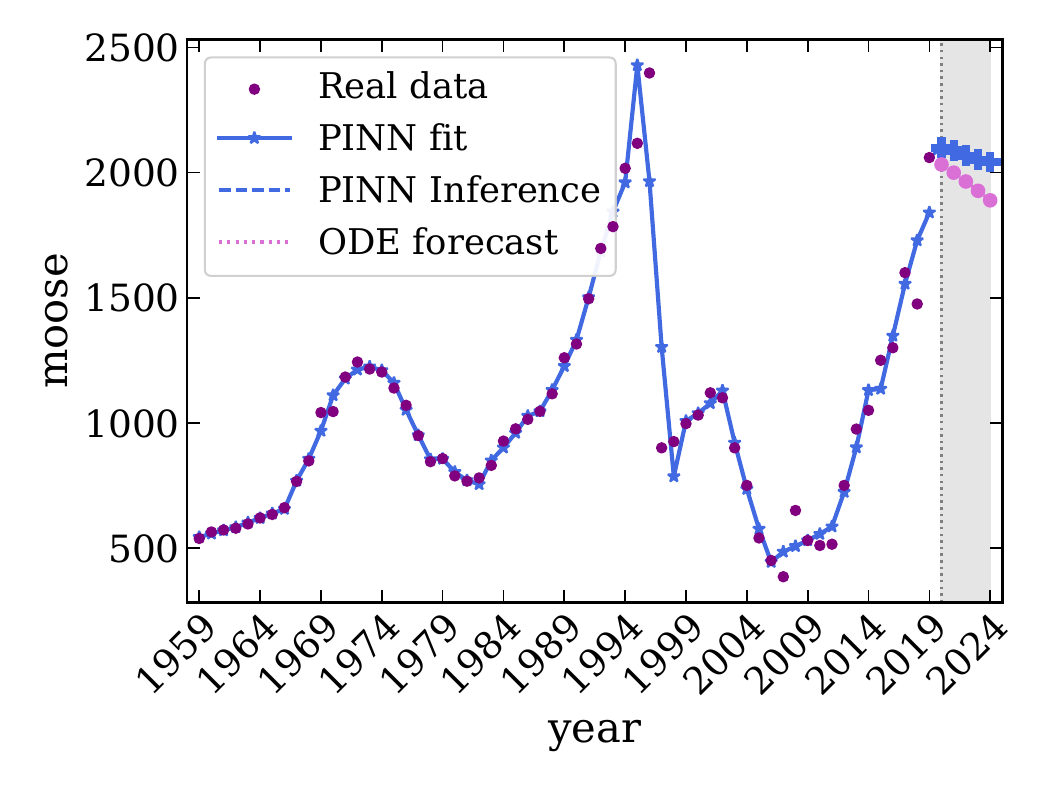}
		\caption{}
		\label{fig-8(a)}
	\end{subfigure}
	\hfill
	\begin{subfigure}{0.49\textwidth}
		\centering
		\includegraphics[width=\linewidth]{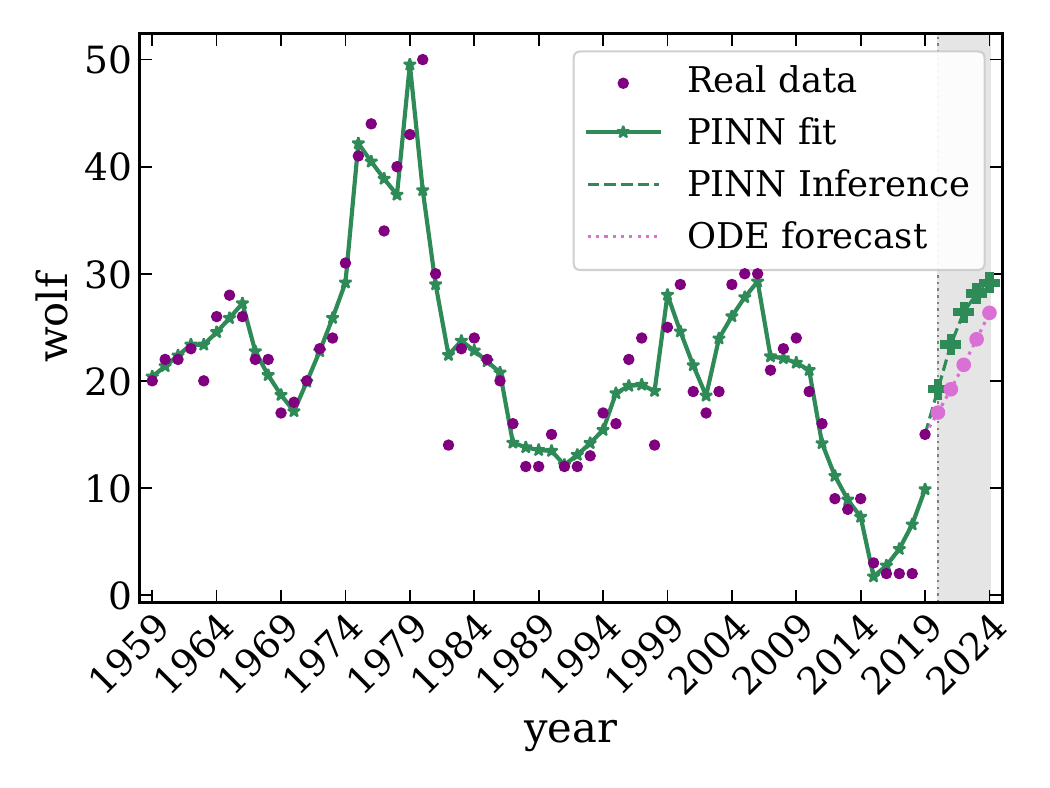}
		\caption{}
		\label{fig-8(b)}
	\end{subfigure}	
	\caption{\textbf{Ratio-dependent model:} Observed and fitted population trends (1959--2019), including predictions (2020--2024). Panel \textbf{(a)} shows the moose population dynamics, while panel \textbf{(b)} shows the wolf population dynamics.}
	\label{fig-8}
\end{figure}

The neural network for the ratio-dependent model was trained independently on each segment using the Adam optimizer, with $15000$ iterations performed for every segment. Since the entire time domain was partitioned into $15$ segments, the total number of training iterations reached $225000$.
Fig.~\ref{fig-8} presents the observed and model fitted population dynamics of moose and wolves on Isle Royale over a six-decade period $(1959-2019)$, along with model generated predictions extending through $2024$. Fig.~\ref{fig-8(a)} illustrates the moose population trajectory, which exhibits strong fluctuations over some years. The fitted model closely tracks these observed oscillations. Predictions for $2020-2024$ suggest a continued declining trend, which follows the actual population trend (see Fig.~\ref{fig-2}). We see the largest error in window $10$ (Table~\ref{table-1}) because the moose population collapsed suddenly due to environmental conditions. 

\begin{figure}[ht!]
	\centering	
	\begin{subfigure}{0.49\textwidth}
		\centering
		\includegraphics[width=\linewidth]{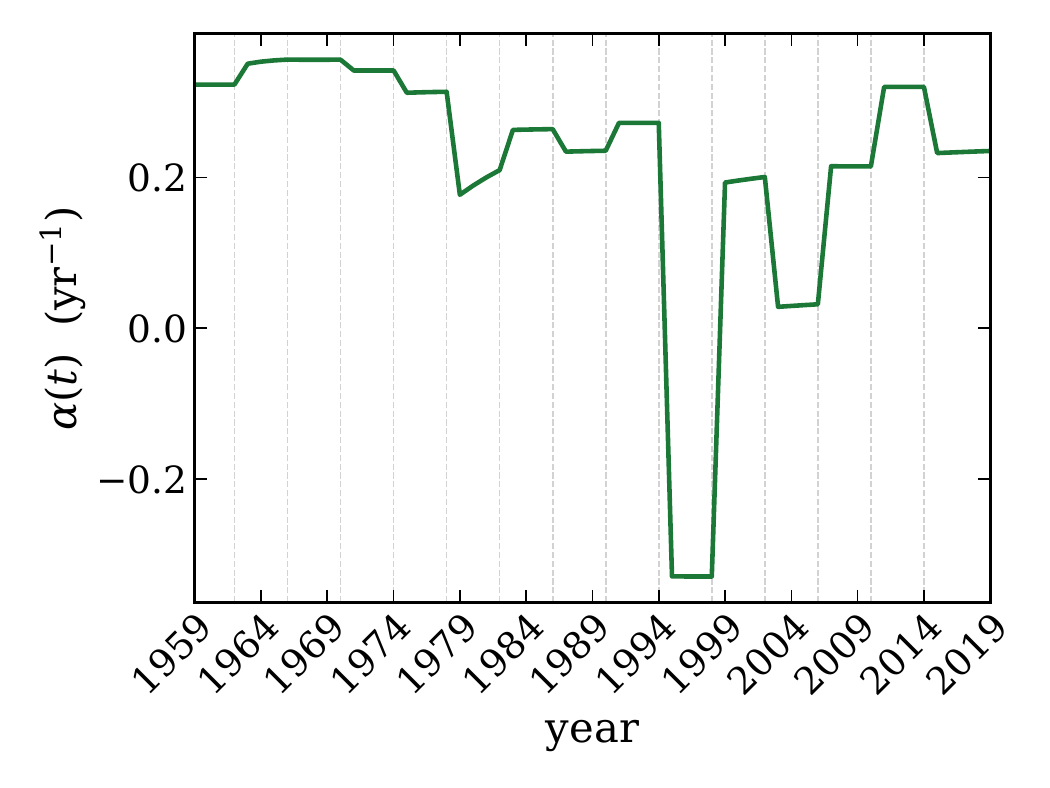}
		\caption{}
		\label{fig-9(a)}
	\end{subfigure}
	\hfill
	\begin{subfigure}{0.49\textwidth}
		\centering
		\includegraphics[width=\linewidth]{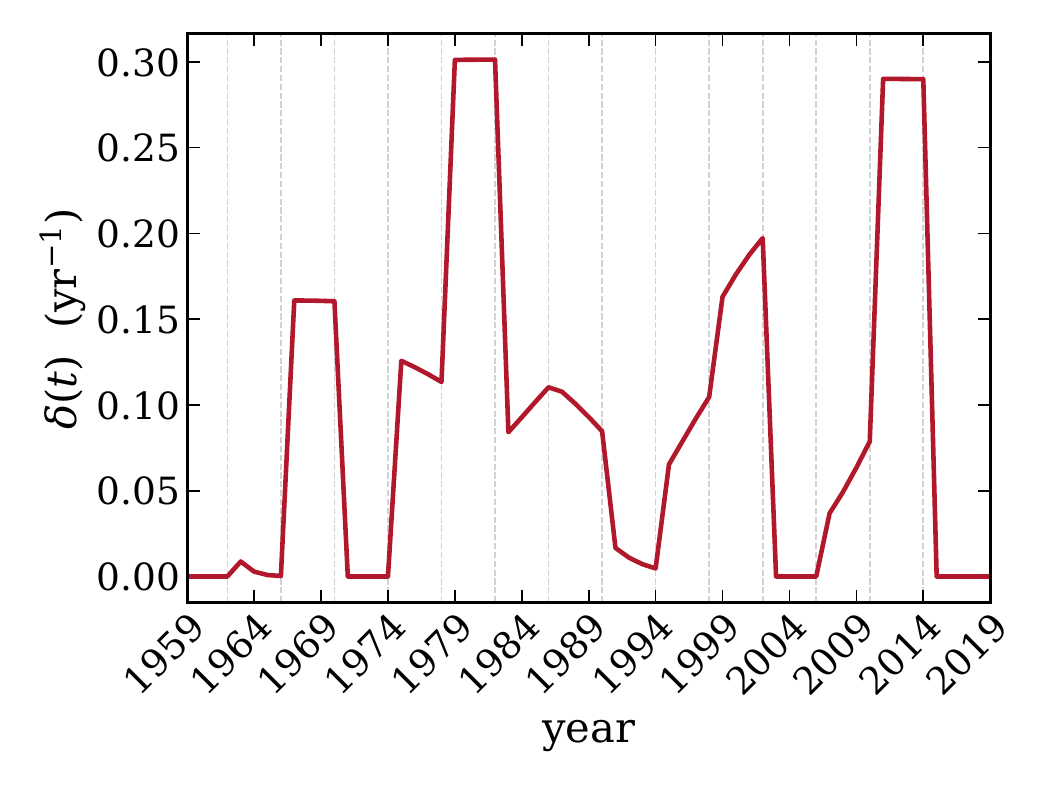}
		\caption{}
		\label{fig-9(b)}
	\end{subfigure}
	\caption{Estimated time dependent parameters for \textbf{ratio-dependent model}: \textbf{(a)} $\alpha(t)$ and \textbf{(b)} $\delta(t)$.}
	\label{fig-9}
\end{figure}

Fig.~\ref{fig-8(b)} depicts the corresponding wolf population dynamics over the same period. The PINN model is able to interpolate the real training data and giving good prediction outside the training set. 

\begin{figure}[ht!]
	\centering	
	\begin{subfigure}{0.32\textwidth}
		\centering
		\includegraphics[width=\linewidth]{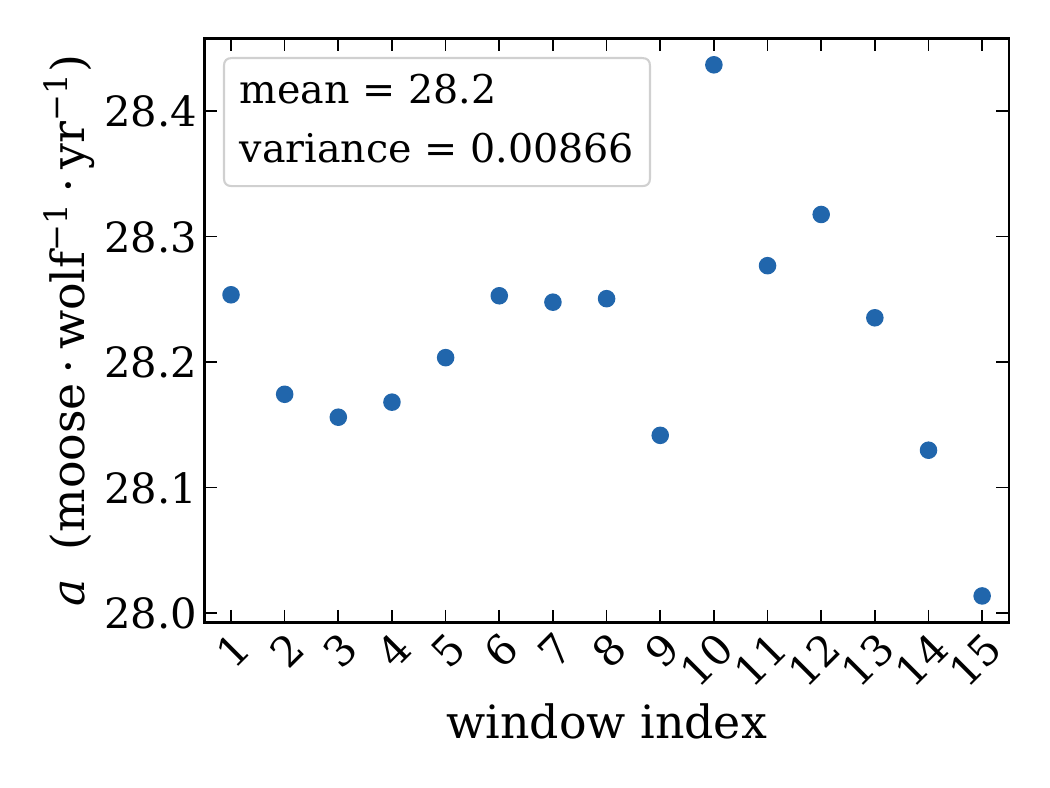}
		\caption{}
		\label{fig-10(a)}
	\end{subfigure}
	\hfill
	\begin{subfigure}{0.32\textwidth}
		\centering
		\includegraphics[width=\linewidth]{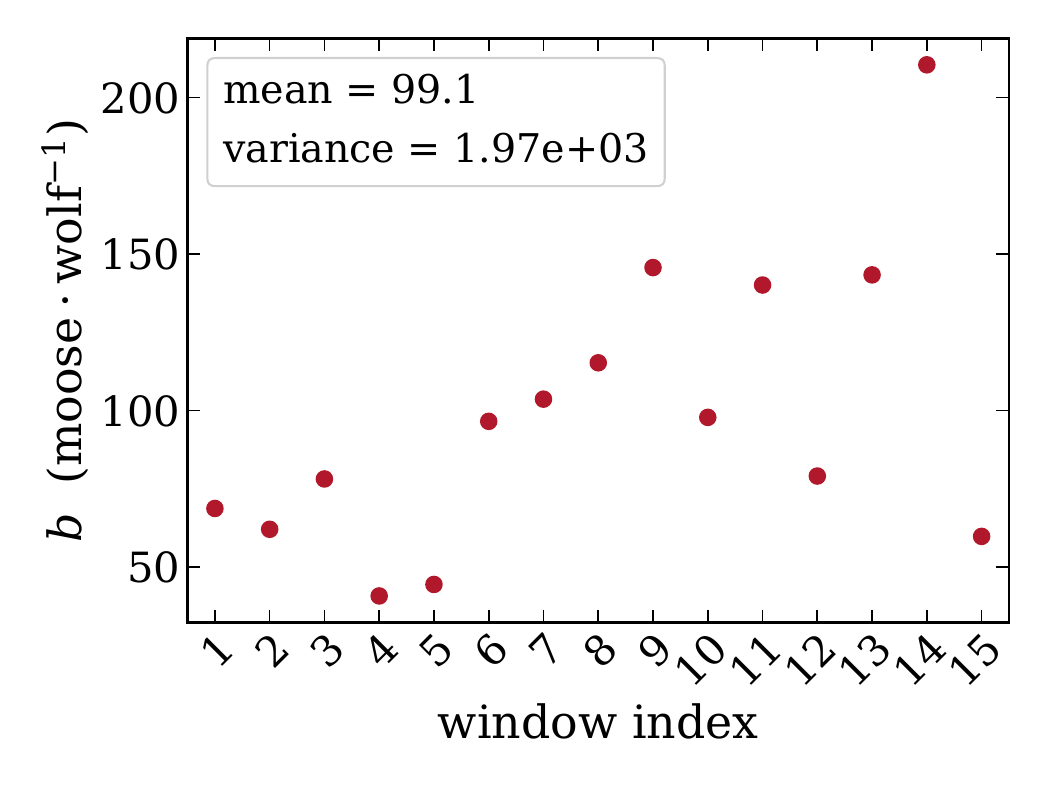}
		\caption{}
		\label{fig-10(b)}
	\end{subfigure}
	\hfill	
	\begin{subfigure}{0.32\textwidth}
		\centering
		\includegraphics[width=\linewidth]{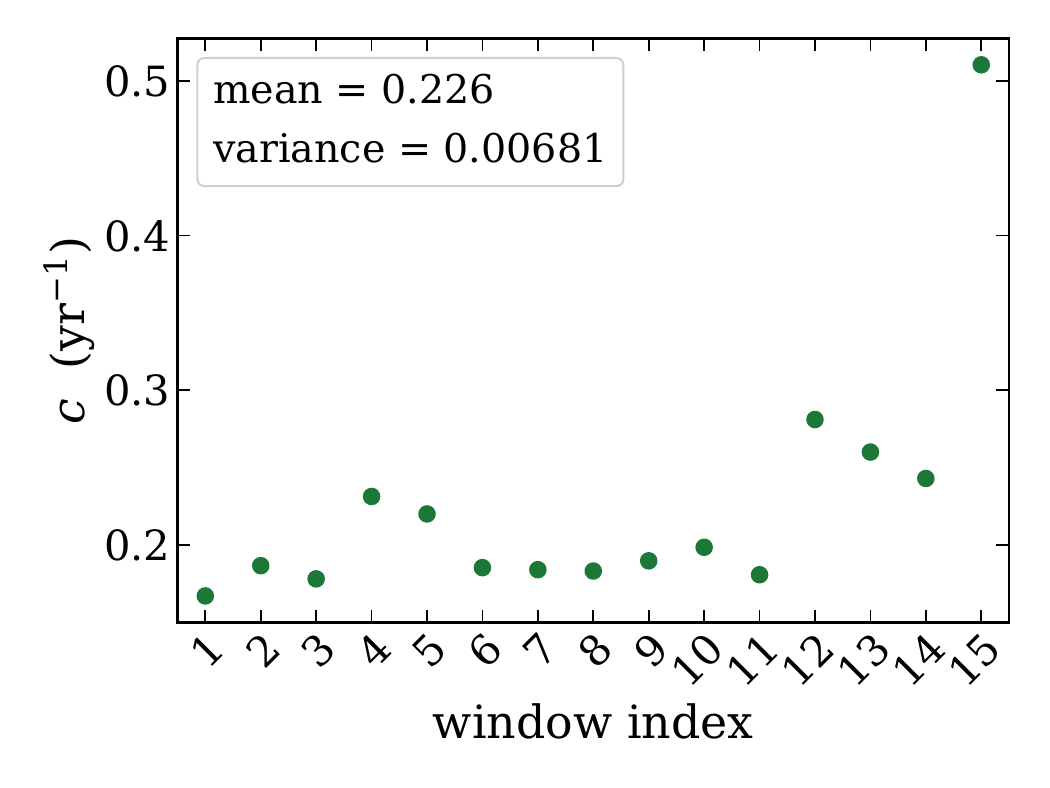}
		\caption{}
		\label{fig-10(c)}
	\end{subfigure}
	\caption{Estimated constant parameters for \textbf{ratio-dependent model}--a, b, and c across all window with their respective mean and variance.}
	\label{fig-10}
\end{figure}

\begin{figure}[ht!]
	\centering	
	\begin{subfigure}{0.49\textwidth}
		\centering
		\includegraphics[width=\linewidth]{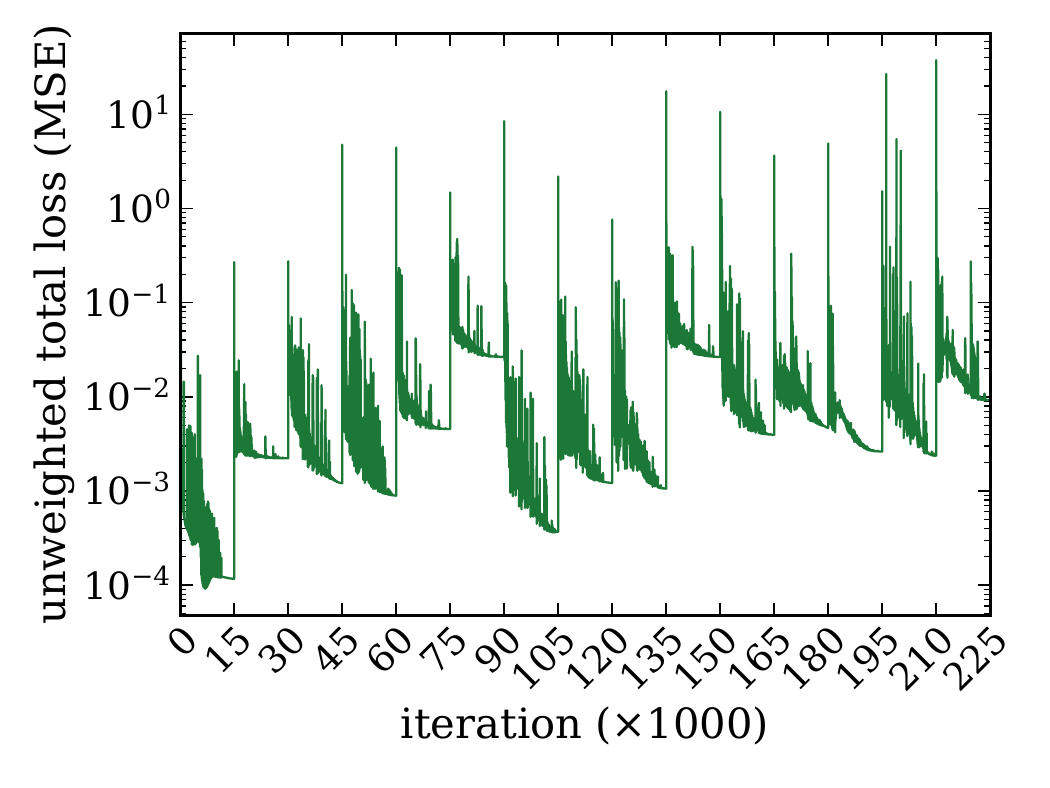}
		\caption{}
		\label{fig-11(a)}
	\end{subfigure}
	\hfill
	\begin{subfigure}{0.49\textwidth}
		\centering
		\includegraphics[width=\linewidth]{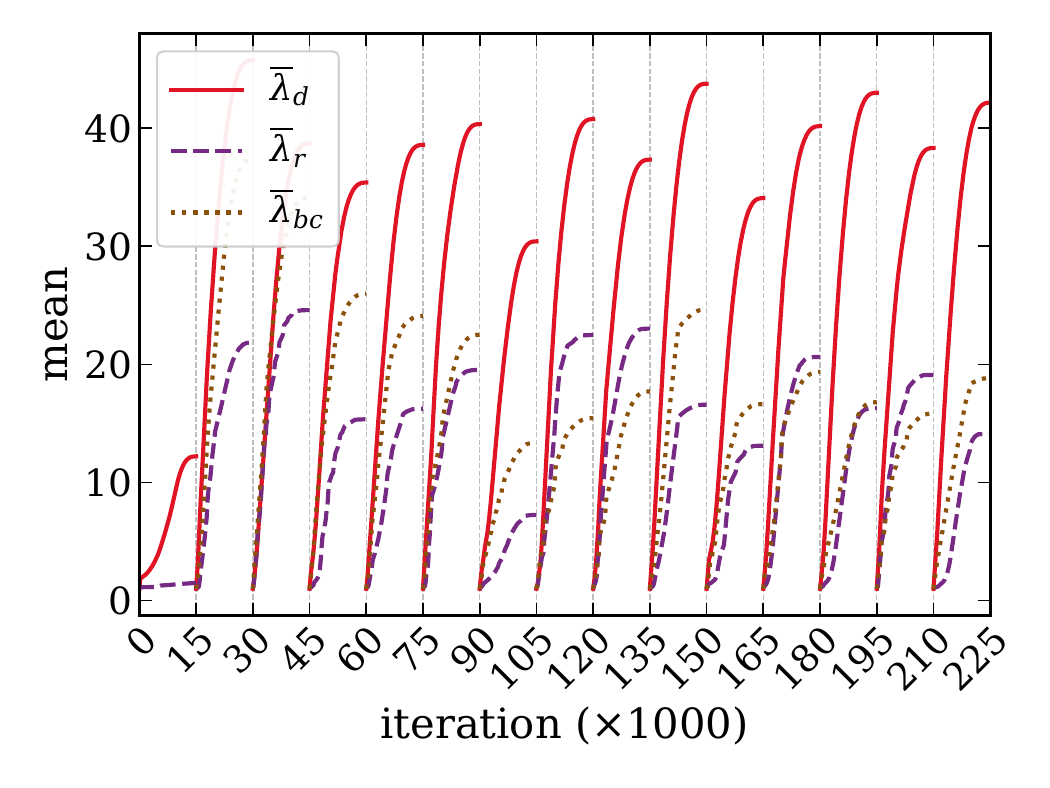}
		\caption{}
		\label{fig-11(b)}
	\end{subfigure}
	\caption{Convergence history of the sequential sweep for \textbf{ratio-dependent model}, plotted against the cumulative
		iteration count over all $15$ segments ($15000$ Adam iterations each, $225000$ in
		total).
		\textbf{(a)} Unweighted total loss, $\mathcal{L}^{\mathrm{unw}}=\mathcal{L}^{\mathrm{unw}}_{d}+
		\mathcal{L}^{\mathrm{unw}}_{r}+\mathcal{L}^{\mathrm{unw}}_{bc}$, obtained by setting
		$\lambda\equiv1$, on a logarithmic scale. Within each segment the error falls by two to
		three orders of magnitude, from order $10^{-1}$ at the transfer to order $10^{-4}$ to
		$10^{-3}$ at convergence; the transient spike at each boundary is the cost of moving the
		inherited network onto a new interval, and its decreasing height in later segments shows
		that transfer learning progressively improves the starting point.
		\textbf{(b)} The average weight given to each of the three parts of the loss: the data term, the
		equation residual term, and the backward-compatible term. These weights are not fixed
		numbers. During training the weight of any point that is still badly fitted is pushed
		upward, so every weight climbs while a segment is being solved and then drops back to one
		when the next segment begins. The data weights climb highest, to about 40--50, followed by
		the backward-compatible weights and then the residual weights.}
	\label{fig-11}
\end{figure}

The model (\ref{eq-14}) does not consider the disease outbreak in wolves, which was happened in 1981. Also, we see the largest $\ell_2$ relative error (Table \ref{table-1}) which occurs in $15^{th}$ window by which we are making predictions. The main reason for this can be the events, like in $2018-2019$ wolves have been relocated to Isle Royale from other places, which is also not considered by the model (\ref{eq-14}). Despite these shortcomings in the model, the PINN model and ODE model is able to follow the trend of the wolf population from $2020-2024$.

In Fig.~\ref{fig-9}, we see the evolution of the intrinsic growth rate of prey and the natural death rate of predators with time. The growth rate varies from $-0.3293$ to $0.3564$. The ups and downs in growth and death rates influence the model predictions. The identified parameters $a$, $b$, and $c$ obtained via the PINN-based framework across $15$ independent time windows are illustrated in Fig.~\ref{fig-10}. All three parameters consistently oscillate around their respective mean values $28.2$, $99.1$, and $0.226$, with variances $0.008$, $1.97\times10^{3}$, and $0.006$, respectively. Fig.~\ref{fig-11} presents the convergence of the bc-PINN and evolution of the mean weight given to the each loss component.

\subsection{Comparison of results}

We have implemented a self-adaptive bc-PINN with transfer learning framework to model systems (\ref{eq-9}) and (\ref{eq-14}). Both models show their ability in terms of predicting the population dynamics of moose and wolves. From the results, we see that the ratio-dependent formulation is more capable of predicting the population of moose than type-II (Figs.~\ref{fig-2}, \ref{fig-4(a)}, and \ref{fig-8(a)}). In the case of wolves, both models are following the upward trend (Figs.~\ref{fig-2}, \ref{fig-4(b)}, and \ref{fig-8(b)}). Also Table~\ref{table-1} shows the $\ell_2$ relative errors and $R^{2}$ score for both the models. From this, we can also see that the overall error is lower for the ratio-dependent model than for type-II; here also, the ratio-dependent model performs better. The moose population growth rate from 1959--1998 on Isle Royale is also reported in Fig.~3 of \citet{vucetich2004influence} which approximately matches with our estimated growth rates. Segment by segment the ratio-dependent form is the more accurate in $10$ of $15$ segments for moose and $11$ of $15$ for wolves. Both responses share the same difficult segments: segment $6$, which spans the collapse of wolves due to disease in $1981$, segment $10$ ($1995$--$1998$), which spans the collapse of the moose, and segment $15$ ($2015$--$2019$), which spans the subsequent reintroduction of the wolves (see Table~\ref{table-3}). 

\begin{table}[ht!]
	\centering
	\caption{Accuracy of the bc-PINN reconstruction of the Isle Royale moose and wolves record,
		compared between the two functional responses. For each of the $15$ sequential time
		segments, and for the complete record, the table reports the relative $\ell_2$ error,
		$\varepsilon_{\ell_2}=\lVert \hat{v}-v\rVert_2/\lVert v\rVert_2$ and the coefficient of
		determination $R^2=1-\sum_k(\hat v_k-v_k)^2/\sum_k(v_k-\bar v)^2$, computed separately for
		the moose population ($m$) and the wolf population ($w$). Segments $1$ to $14$ contain four
		annual observations each and segment $15$ contains five; the full-record statistics use all $61$ observations. Over the complete record the two responses are close, the ratio-dependent form giving a
		$3.9\%$ smaller relative error for moose ($0.0949$ against $0.0987$) and an identical wolf
		error to four decimal places. }
	\label{table-1}
	\hfill
	\begin{minipage}[t]{0.54\textwidth}
		\centering
		(a) Holling type-II response\\[4pt]
		\small
		\begin{tabular}{r l r r r r}
			\hline
			\multicolumn{2}{c}{} & \multicolumn{2}{c}{Moose ($m$)} & \multicolumn{2}{c}{Wolf ($w$)} \\
			\cmidrule(lr){3-4}\cmidrule(lr){5-6}
			Win. & Years & $\varepsilon_{\ell_2}$ & $R^2$ & $\varepsilon_{\ell_2}$ & $R^2$ \\
			\hline
			1  & 1959--1962 & 0.0086 & 0.9027 & 0.0247 & 0.7558 \\
			2  & 1963--1966 & 0.0080 & 0.9547 & 0.0887 & 0.4451 \\
			3  & 1967--1970 & 0.0541 & 0.8274 & 0.0643 & 0.6847 \\
			4  & 1971--1974 & 0.0126 & 0.5076 & 0.0577 & 0.8736 \\
			5  & 1975--1978 & 0.0228 & 0.9587 & 0.0852 & 0.1227 \\
			6  & 1979--1982 & 0.0122 & 0.9214 & 0.2190 & 0.6531 \\
			7  & 1983--1986 & 0.0215 & 0.9141 & 0.0363 & 0.7000 \\
			8  & 1987--1990 & 0.0164 & 0.9674 & 0.1279 & 0.0135 \\
			9  & 1991--1994 & 0.0225 & 0.9551 & 0.0923 & 0.6260 \\
			10 & 1995--1998 & 0.2079 & 0.7219 & 0.2067 & 0.0502 \\
			11 & 1999--2002 & 0.0291 & 0.6171 & 0.1417 & 0.5334 \\
			12 & 2003--2006 & 0.0451 & 0.9695 & 0.0954 & 0.6824 \\
			13 & 2007--2010 & 0.1689 & 0.1023 & 0.0645 & 0.4617 \\
			14 & 2011--2014 & 0.0771 & 0.9020 & 0.1672 & 0.6712 \\
			15 & 2015--2019 & 0.1051 & 0.6784 & 0.4593 & 0.6032 \\
			\hline
			\multicolumn{2}{l}{\textbf{Full record}} & \textbf{0.0987} & \textbf{0.9357}
			& \textbf{0.1316} & \textbf{0.9045} \\
			\hline
		\end{tabular}
	\end{minipage}
	\hfill
	\begin{minipage}[t]{0.444\textwidth}
		\centering
		(b) Ratio-dependent response\\[4pt]
		\small
		\begin{tabular}{r r r r}
			\hline
			\multicolumn{2}{c}{Moose ($m$)} & \multicolumn{2}{c}{Wolf ($w$)} \\
			\cmidrule(lr){1-2}\cmidrule(lr){3-4}
			$\varepsilon_{\ell_2}$ & $R^2$ & $\varepsilon_{\ell_2}$ & $R^2$ \\
			\hline
		     0.0084 & 0.9063 & 0.0210 & 0.8243 \\
			 0.0067 & 0.9676 & 0.0880 & 0.4545 \\
			 0.0526 & 0.8371 & 0.0630 & 0.6980 \\
			 0.0137 & 0.4145 & 0.0527 & 0.8946 \\
			 0.0144 & 0.9834 & 0.0835 & 0.1580 \\
			 0.0187 & 0.8168 & 0.2199 & 0.6502 \\
			 0.0207 & 0.9208 & 0.0359 & 0.7065 \\
			 0.0169 & 0.9655 & 0.1203 & 0.1271 \\
			 0.0236 & 0.9505 & 0.0828 & 0.6990 \\
			 0.1980 & 0.7478 & 0.1973 & 0.1348 \\
			 0.0246 & 0.7270 & 0.1324 & 0.5926 \\
			 0.0326 & 0.9840 & 0.1136 & 0.5493 \\
			 0.1711 & 0.0787 & 0.0787 & 0.1996 \\
			 0.0780 & 0.8997 & 0.1549 & 0.7180 \\
			 0.1031 & 0.6905 & 0.4720 & 0.5811 \\
			\hline
			\textbf{0.0949} & \textbf{0.9405} & \textbf{0.1316} & \textbf{0.9045} \\
			\hline
		\end{tabular}
	\end{minipage}
\end{table}

\begin{table}[ht!]
	\centering
	\caption{\textbf{Comparison of estimated parameters across studies.} The year $1971$ marks the first availability of annual estimates of wolf kill rates. Our parameters were estimated using prey--predator mathematical models, whereas the other two studies focused only on functional responses. However, the parameters have the same ecological meaning in all three studies. The bounds for parameter $a$ in our estimation were selected based on the ranges reported in the previous studies. While the previous studies reported ranges of parameter values, we report the mean values estimated using the bc-PINN model across $15$ segments.}
	\label{table-2}
	\begin{tabular}{llcccc}
		\toprule
		\multirow{2}{*}{Method} & \multirow{2}{*}{Study period}
		& \multicolumn{2}{c}{Ratio-dependent} & \multicolumn{2}{c}{Holling type-II} \\
		\cmidrule(lr){3-4} \cmidrule(lr){5-6}
		&            & $a$ & $b$   & $a$ & $b$ \\
		\midrule
  	\textbf{1.} PINN method        & 1959--2019 & 28.217 & 99.061 & 28.175  & 3332.692  \\
	\textbf{2.} \citet{vucetich2002effect}  & 1971--2001 & [18.62,~23.26] & [31.41,~56.38]   & [21.81,~41.37]  & 										[1001.70,~3011.52]  \\
   \textbf{3.} \citet{jost2005wolves}      & 1971--1998 & [15.58,~31.46] & [47.44,~179.58]   & --- & --- \\
		\bottomrule
	\end{tabular}
\end{table}

\textbf{Remark:} The real time units in Method 2 were months, whereas those in Method 3 were days; we converted them to years for comparison.
 
\section{Discussion and conclusion} \label{sec-6}
The study of the moose-wolf system of Isle Royale had been started since $1958$, and many researchers, including ecologists, statisticians, applied mathematicians, and park staff, contributed significantly through these years. Their collaborative efforts have led to a more profound understanding of the complex interactions between these species and the overall ecosystem. As a result, the findings from this long-term research continue to inform conservation strategies and wildlife management practices. The role of Michigan Technological University researchers in the moose-wolf project has shaped this closed habitat into a research environment for others. The collective efforts made by the researchers have found many interesting results from ecological and behavioral perspectives, but there is a lack of models that can predict the dynamics in the future. To fill this gap, we have combined the data and previous studies on Isle Royale, such as functional response studies and types of logistic growth studies, with trending deep learning algorithms to create a model that can predict the future population of moose and wolves. Based on the best-performing functional responses, we build two economical models on the basis of observations. We have only population data for moose and wolves; therefore, we are unable to increase the number of parameters. Since we are solving an inverse problem, increasing the number of parameters without additional data poses significant challenges. Parameter estimation with PINN is inherently an ill-posed inverse problem, and this becomes increasingly difficult as more parameters are added without additional constraints, regularization, or prior information. Due to data limitations, we tried the simplest model with fewer parameters.The limitations of both models highlight the need for additional information to improve the accuracy of future population predictions and parameter estimation, thereby providing more reliable insights to support informed decision-making by park administrators.

\noindent \textbf{Limitations of the considered models:}
\begin{itemize}
	\item Both models assume a well-mixed, closed system with no immigration, emigration, ignoring the spatially structured, pack-organized predation and age-selective prey vulnerability that characterize the Isle Royale system.
	
	\item Both models treat predator dynamics as purely endogenous; however, 19 wolves were released at Isle Royale in 2018--2019 (Table {\ref{table-3}}) by the National Park Service, constituting an abrupt, externally imposed discontinuity in wolf abundance that no smooth ODE formulation can reproduce. We attribute the largest prediction error (Table~\ref{table-1}) in the final observation window to this human intervention.
	
	\item Both models capture wolf mortality through a function $\delta(t)$, yet canine parvovirus, unintentionally introduced by humans, was largely responsible for the wolf population crashing from 50 to 14 individuals, a rapid, pathogen-driven collapse incoherent with any continuous mortality term. Embedding an epidemiological structure within the predator equation would substantially improve predictive fidelity during outbreak periods.
	
	\item A more faithful representation of the Isle Royale moose--wolf dynamics would require additional data, which are currently unavailable to us and therefore beyond the scope of the present study.
\end{itemize}

Despite these above-mentioned limitations, our models are able to interpolate and extrapolate. With our results, we see that ratio-dependent formulations have the better capability of extrapolation than type-II. Our results show that the ratio-dependent model forecast the population trend of moose and wolf very well, and it can be verified using graphical data (Fig.~\ref{fig-2}).

\begin{table}[ht!]
	\centering
	\caption{Key events and moose-wolf population changes at Isle Royale National Park (\href{https://www.isleroyalewolf.org}{https://www.isleroyalewolf.org}) from 1959--2024.}
	\label{table-3}
	\begin{tabular}{llcc}
		\hline
		\textbf{Year} & \textbf{Key Event} & \textbf{Wolves trend} & \textbf{Moose trend} \\
		\hline
		1959--1968  & Both animals lived in balance           & Stable      & Stable         \\
		$\sim$1973  & Good conditions for moose               & Moderate    & Nearly doubled \\
		$\sim$1980  & Hard winters; wolves hunted more        & Increased   & Fell by half   \\
		1981--1982  & Disease (Canine parvovirus) hit wolves  & Crashed     & Rose           \\
		1997        & Winter ticks, bad winters, less food    & Low         & Collapsed      \\
		$\sim$2018  & Very few wolves left                    & Near zero   & Very high      \\
		2018--2019  & Wolves were brought back by humans      & Started increasing  & Increasing      \\
		After 2019  & After reintroduction of wolves          & Increasing  & Started decreasing \\
		\hline
	\end{tabular}
\end{table}

Overall, the study demonstrates the potential of advanced deep learning frameworks to address complex problems in ecology. The results show that the framework can predict population dynamics outside the range of the training data, although some limitations of the proposed approach are identified. This ability to generalize outside the training domain is especially important as extrapolation remains one of the main challenges in neural network based models. These findings highlight the promise of deep learning as a tool for ecological modeling and provide a foundation for future work to improve model accuracy, interpretability, and robustness. In our implementation in PyTorch, we set a fixed random seed to ensure the reproducibility of results within the same hardware environment. Fully reproducible results are not guaranteed across platforms and results may not be reproducible between CPU and GPU executions even with the same seeds or different seeds (\href{https://docs.pytorch.org/docs/2.12/notes/randomness.html}{https://docs.pytorch.org/docs/2.12/notes/randomness.html}). Thus, different convergence behaviors and final solutions can be obtained when changing the seed -- or hardware, especially for ill-posed inverse problems where the loss landscape is highly non-convex and the optimization is sensitive to initialization and numerical perturbations.

\section*{Authorship contribution statement}

\textbf{Anurag Singh}: Writing – original draft, Writing – review \& editing, Software, Formal analysis, Methodology, Visualization, Conceptualization. \textbf{Nitu Kumari}: Writing – review \& editing, Methodology, Formal analysis, Validation, Visualization, Supervision, Conceptualization. 

\section*{Acknowledgments}
 The research of \textbf{Nitu Kumari} was supported by the Council of Scientific and Industrial Research (CSIR, India), under the research grant $\#$ $25\slash0326\slash23\slash\text{EMR}-\text{II}$. We sincerely thank \textbf{Prof. Khemraj Shukla, Applied Mathematics Division, Brown University, USA} for his continuous guidance and support throughout this work. His valuable insights and assistance in the implementation of physics-informed neural network were instrumental in the successful completion of this study. We also thank \textbf{Prof. John A. Vucetich, Michigan Technological University} for his valuable suggestions.

\section*{Conflict of interest}
The authors declare no conflict of interest.

\section*{Data and code availability}
The data and code will be made publicly available in a permanent repository upon publication.

\bibliographystyle{apalike} 
\bibliography{ref}

\end{document}